\documentclass[12pt]{amsart}

\usepackage{amsmath,amssymb,amsthm}
\usepackage{epsf}
\usepackage{epsfig}
\usepackage[T2A]{fontenc}
\usepackage[cp1251]{inputenc}
\usepackage[english]{babel}

\newcommand{\RR}{\mathbb R}
\newcommand{\CC}{\mathbb C}

\renewcommand{\Re}{\mathop{\rm Re}\nolimits}
\renewcommand{\Im}{\mathop{\rm Im}\nolimits}
\newcommand{\trace}{\mathop{\mathrm{trace}}\nolimits}

\newcommand{\ii}{\mathrm{i}}

\newcommand{\manifold}[1]{\mathcal{#1}}
\newcommand{\M}{\manifold{M}}
\newcommand{\D}{\manifold{D}}

\newcommand{\vect}[1]{\mathrm{#1}} 
\newcommand{\x}{\vect{x}}
\newcommand{\y}{\vect{y}}
\newcommand{\va}{\vect{a}}
\newcommand{\vb}{\vect{b}}
\newcommand{\vc}{\vect{c}}
\newcommand{\vd}{\vect{d}}
\newcommand{\vX}{\vect{X}}
\newcommand{\vY}{\vect{Y}}
\newcommand{\vH}{\vect{H}}
\newcommand{\n}{\vect{n}}

\newtheorem{thm}{Theorem}[section]

\newtheorem{prop}[thm]{Proposition}
\theoremstyle{definition}

\theoremstyle{remark}
\newtheorem{rem}[thm]{Remark}

\newcommand{\ds}{\displaystyle}

\begin{document}

\title [Canonical Weierstrass representations for minimal surfaces in $\RR^4$]
{Canonical Weierstrass representations for minimal surfaces in Euclidean 4-space}%

\author{Georgi Ganchev and Krasimir Kanchev}

\address{Bulgarian Academy of Sciences, Institute of Mathematics and Informatics,
Acad. G. Bonchev Str. bl. 8, 1113 Sofia, Bulgaria}
\email{ganchev@math.bas.bg}%

\address {Department of Mathematics and Informatics, Todor Kableshkov University of Transport,
158 Geo Milev Str., 1574 Sofia, Bulgaria}%
\email{kbkanchev@yahoo.com}%


\begin{abstract}
Minimal surfaces of general type in Euclidean 4-space are characterized with
the conditions that the ellipse of curvature at any point is centered at this point and has
two different principal axes. Any minimal surface of general type locally admits geometrically
determined parameters - canonical parameters. In such parameters the Gauss curvature and
the normal curvature satisfy a system of two natural partial differential equations and
determine the surface up to a motion. For any minimal surface parametrized by canonical parameters we obtain
Weierstrass representations - canonical Weierstrass representations. These Weierstrass formulas
allow us to solve explicitly the system of natural partial differential equations and to establish geometric
correspondence between minimal surfaces of general type, the solutions to the system of
natural equations and pairs of holomorphic functions in the Gauss plane. On the base of these
correspondences we obtain that any minimal surface of general type in Euclidean 4-space
determines locally a pair of two minimal surfaces in Euclidean 3-space and vice versa. Finally
some applications of this phenomenon are given.
\end{abstract}

\maketitle

\thispagestyle{empty}
\section{Introduction}
Analytic methods to study surfaces and their properties are of
essential importance in differential geometry. A classical example
of such an approach is given by the Weierstrass representation for
minimal surfaces, which is one of the most powerful instruments
for constructing new surfaces.

In this paper we study minimal surfaces in the four-dimensional
Euclidean space $\RR^4$. For any surface $\M$ in $\RR^4$ we denote
by $K, \, \varkappa$ and $H$ the Gauss curvature, the normal
curvature and the mean curvature, respectively. These three
invariants satisfy the following inequality \cite{W}
$$K+|\varkappa|\leq ||H||^2.$$
A surface $\M$ is said to be minimal if $H=0$, which means
geometrically that the ellipse of curvature at any point is
centered at this point. Therefore any minimal surface satisfies
the inequality
$$K^2-\varkappa^2\geq 0,$$
which divides the minimal surfaces into two classes:
\begin{itemize}

\item
the class of minimal super-conformal surfaces characterized by
$K^2-\varkappa^2=0$;

\item
the class of minimal surfaces of general type characterized by
$K^2-\varkappa^2>0$.
\end{itemize}

Geometrically, any superconformal surface is characterized by the
condition that its ellipse of curvature is a circle. Minimal
superconformal surfaces in $\RR^4$ were described geometrically in
\cite{DT} (see also \cite{M}).

Here we consider minimal surfaces in $\RR^4$ of general type. At
any point of such a surface the ellipse of curvature has two
different principal axes.

Next we describe our scheme of investigation.

The leading idea is to study surfaces in $\RR^3$ or in $\RR^4$
with respect to special geometrically determined parameters -
\emph{canonical parameters} \cite {GM}. With respect to such
parameters all coefficients of the first and the second
fundamental form are expressed by the invariants of the surface.

Any minimal surface in $\RR^4$ of general type admits special isothermal
parameters - canonical parameters (of the first type or of the
second type) (cf \cite{Itoh}). To endow locally the minimal
surface under consideration with these parameters means that the
tangent to any parametric line is transformed by the second
fundamental tensor in a normal, which is collinear to a principal
axis of the ellipse of curvature. Further we obtain Weierstrass
representation formulas with respect to canonical parameters which
describe locally all minimal surfaces in terms of two holomorphic
functions. Introducing canonical parameters on a minimal surface,
one obtains the system of natural PDE's of minimal surfaces and a
Bonnet type fundamental theorem for minimal surfaces of general
type \cite{TG}. The canonical Weierstrass formulas allow us to
obtain explicitly the solutions to the system of natural PDE's of
minimal surfaces \cite{GK1}.

We consider the set $\mathbf{MS_4}$ of equivalence classes of
minimal surfaces containing a fixed point, the set
$\mathbf{SNE_4}$ of equivalence classes of solutions to the system
of natural PDE's and the set $\mathbf{H^2}$ of equivalent pairs of
holomorphic functions. Our main result is that any two of these
sets $\{\mathbf{MS_4}, \; \mathbf{SNE_4}, \; \mathbf{H^2}\}$ are
in a natural one-to-one correspondence.

This result leads to a natural correspondence between the minimal
surfaces in $\RR^4$ and pairs of minimal surfaces in $\RR^3$,
which is a base of a systematical study of minimal surfaces in
$\RR^4$ having in mind the well developed theory of minimal
surfaces in $\RR^3$.

\section{Preliminaries}
Let $\M$ be a two-dimensional Riemannian manifold and $\x : \M \to
\RR^n$ be an isometric immersion of $\M$ into $\RR^n$. Then we say
that $(\M ,\x)$ (or $\M$) is a regular surface in $\RR^n$. If $\x
: (u,v) \to \x(u, v) \in \RR^n; \, (u,v)\in \D \subset \RR^2$ is a
parametrization of $\M$, then the coefficients of the first
fundamental form are $E=\x_u^2$, $F=\x_u\cdot\x_v$ and $G=\x_v^2$.
Without loss of generality we can assume that the parameters
$(u,v)$ are isothermal local coordinates, i.e. $E=G$ and $F=0$.

In addition to the real coordinates $(u,v)$ we also consider the
complex coordinate $t=u+\ii v$, identifying the coordinate plane
$\RR^2$ with the complex plane $\CC$. Thus all functions defined
on the surface can be considered as functions of the complex
variable $t$.

We denote by $T_p(\M)$ the tangent plane of $\M$ at a point $p \in
\M$, which is identified with the corresponding plane in $\RR^n$.
The normal space $N_p(\M)$ at $p$ is the normal complement of
$T_p(\M)$ in $\RR^n$. Using the standard imbedding of $\RR^n$ into
$\CC^n$ we consider the complexified tangent space $T_{p,C}(\M)$
to $\M$ at the point $p$ as a subspace of $\CC^n$, which is the
linear span of $T_p(\M)$ in $\CC^n$. In a similar way the
complexified normal space $N_{p,C}(\M)$ to $\M$ is identified with
the corresponding subspace of $\CC^n$, which is the linear span of
$N_p(\M)$ in $\CC^n$.

If $\va$ and $\vb$ are two vectors in $\CC^n$, then $\va\cdot \vb$
(or $\va\vb$) denotes the  bilinear dot product
\[\va\cdot \vb=a_1b_1+a_2b_2+\cdots +a_nb_n\; \]
and the dot product of the vector $\va$ with itself is
\[\va^2=\va\cdot \va=a_1^2+a_2^2+\cdots +a_n^2\; .\]
The Hermitian dot product of $\va$ and $\vb$ is given by the
formula
\[\va\cdot\bar \vb = a_1\bar b_1+a_2\bar b_2+\cdots +a_n\bar b_n\;\]
and the norm of the vector $\va$ with respect to the Hermitian dot
product is
\[\|\va\|^2=\va\cdot\bar \va=|a_1|^2+|a_2|^2+\cdots +|a_n|^2\; .\]

Since $T_{p,C}(\M)$ and $N_{p,C}(\M)$ are generated by the real
spaces $T_p(\M)$ and $N_p(\M)$ respectively, they are closed under
the complex conjugation. They are mutually orthogonal with respect
to the bilinear or Hermitian dot product. Therefore we have the
following orthogonal decomposition:
\[\CC^n = T_{p,C}(\M) \oplus N_{p,C}(\M)\; .\]
For a given vector $\va$ in $\CC^n$ $\va^\top$ and $\va^\bot$
denote the orthogonal projections of $\va$ into $T_{p,C}(\M)$ and
$N_{p,C}(\M)$, respectively. For any vector we have:
\[\va=\va^\top + \va^\bot\; .\]
This decomposition is valid with respect to both dot products in
$\CC^n$.

The second fundamental form $\sigma$ of $\M$, is given by:
\[\sigma (\vX,\vY) = (\nabla_\vX \vY)^\bot \; , \]
where $\vX, \, \vY\in T(\M)$, and $\nabla$ is the canonical linear
connection in $\RR^n$.

Let $\vX_1$ and $\vX_2$ be the unit coordinate vector fields on
$\M$ of the same direction as $\x_u$ and $\x_v$, respectively,
i.e.
\begin{equation*}
\vX_1=\frac{\x_u}{\|\x_u\|}=\frac{\x_u}{\sqrt E};\quad
\vX_2=\frac{\x_v}{\|\x_v\|}=\frac{\x_v}{\sqrt G}=\frac{\x_v}{\sqrt
E}\; .
\end{equation*}

As usual, $\vH$ will denote the mean curvature vector field of
$\M$:
\[\vH = \frac{1}{2}\trace\sigma = \frac{1}{2}(\sigma(\vX_1,\vX_1)+\sigma(\vX_2,\vX_2)) \; . \]

Any regular surface with zero mean curvature vector field is said
to be a minimal surface.

\section{The function $\Phi (t)$}

Let $\M: \, \x = \x(u, v); \; (u,v)\in \D \subset \RR^2$ be a
regular surface in $\RR^n$. The complex vector-valued function
$\Phi(t)$ is defined by the equality:
\begin{equation}\label{Phi_def}
\Phi(t)=2\frac{\partial\x}{\partial t}=\x_u-\ii\x_v \;.
\end{equation}

The defining equality implies immediately that
$$\Phi^2=0 \  \Leftrightarrow\  \begin{array}{l} \x_u^2-\x_v^2=0\\  \x_u \x_v=0 \end{array} \
\Leftrightarrow \  \begin{array}{l} E=\x_u^2=\x_v^2=G\\ F=0
\end{array}.$$ Hence, the parameters $(u,v)$ are isothermal if and
only if
\begin{equation}\label{Phi2}
\Phi^2=0\; .
\end{equation}

The norm of $\Phi$ satisfies the following equalities:
\[
\|\Phi\|^2=\Phi\bar{\Phi}=\x_u^2+\x_v^2=E+G=2E=2G.
\]
The coefficients of the first fundamental form of $\M$ are given
in terms of $\Phi$ as follows:
\begin{equation}\label{EG}
E=G=\frac{1}{2}\|\Phi\|^2;\quad F=0 \; .
\end{equation}
Denoting by $\mathbf{I}$ the first fundamental form, then we
have:
\begin{equation}\label{Idt}
\mathbf{I}=\frac{1}{2}\|\Phi\|^2 (du^2 +
dv^2)=\frac{1}{2}\|\Phi\|^2|dt|^2 \; .
\end{equation}
It follows that the function $\Phi$ satisfies the condition:
\begin{equation}\label{modPhi2}
\|\Phi\|^2\neq 0\; .
\end{equation}

Differentiating \eqref{Phi_def} and taking into account the
equality $\frac{\partial}{\partial\bar t} \frac{\partial}{\partial
t} = \frac{1}{4}\Delta$, we find:
\begin{equation}\label{dPhi_dbt}
\frac{\partial\Phi}{\partial\bar t}= \frac{\partial}{\partial\bar
t}\, \left(2\, \frac{\partial \x}{\partial t} \right)=
\frac{1}{2}\Delta \x \; ,
\end{equation}
where $\Delta$ is the Laplace operator.

The last formula implies that $\ds\frac{\partial\Phi}{\partial\bar
t}$ is a real vector-valued function, i.e.
\begin{equation}\label{dPhi_dbt=dbPhi_dt}
\frac{\partial\Phi}{\partial\bar
t}=\frac{\partial\bar\Phi}{\partial t} \; .
\end{equation}

 Thus, any function $\Phi$ given by \eqref{Phi_def} satisfies the conditions:
\eqref{Phi2}, \eqref{modPhi2} and \eqref{dPhi_dbt=dbPhi_dt}.

Conversely, any function $\Phi$ satisfying these three conditions
determines locally the surface up to a translation.

The last assertion follows immediately from the fact that the
condition
\begin{equation}\label{Phi_cond}
\frac{\partial\Phi}{\partial\bar
t}=\frac{\partial\bar\Phi}{\partial t} \;
\end{equation}
is the integrability condition for the system
\begin{equation}\label{xuxv}
\begin{array}{ll}
\x_u=\ \ \,\Re (\Phi)\\[4mm]
\x_v=-\Im (\Phi).
\end{array}
\end{equation}

Further we express the components of the second fundamental form
$\sigma$ by the function $\Phi$.

Taking into account \eqref{dPhi_dbt}, we find:
\begin{equation*}
\left(\frac{\partial\Phi}{\partial\bar t}\right)^\bot=
\left(\frac{1}{2}\Delta \x \right)^\bot= \frac{1}{2}(\x_{uu}^\bot
+ \x_{vv}^\bot)= \frac{1}{2}({\nabla_{\x_u}^\bot \x_u} +
{\nabla_{\x_v}^\bot \x_v })= \frac{1}{2}(\sigma (\x_u,\x_u)+\sigma
(\x_v,\x_v)) .
\end{equation*}

Differentiating \eqref{Phi_def} with respect to $t$ we get:
\begin{equation}\label{dPhi_dt}
\frac{\partial\Phi}{\partial t}= \frac{1}{2}(\x_{uu} - \x_{vv}) -
\ii \x_{uv} .
\end{equation}
Therefore
\begin{equation}\label{dPhi_dt_bot}
\left(\frac{\partial\Phi}{\partial t}\right)^\bot=
\frac{1}{2}(\sigma (\x_u,\x_u)-\sigma (\x_v,\x_v)) - \ii \sigma
(\x_u,\x_v) .
\end{equation}
Hence
\begin{equation}
\begin{array}{l}\label{sigma_uu_vv_uv}
\ds\sigma (\x_u,\x_u) = \Re\left(\frac{\partial\Phi}{\partial\bar
t}\right)^\bot +
                        \Re\left(\frac{\partial\Phi}{\partial     t}\right)^\bot;\\
\ds\sigma (\x_v,\x_v) = \Re\left(\frac{\partial\Phi}{\partial\bar
t}\right)^\bot -
                        \Re\left(\frac{\partial\Phi}{\partial     t}\right)^\bot;\\
\ds\sigma (\x_u,\x_v) = -\Im\left(\frac{\partial\Phi}{\partial
t}\right)^\bot.
\end{array}
\end{equation}

Next we find how the function $\Phi$ is being transformed under a
change of the isothermal coordinates and under a motion of the
surface $(\M,\x)$ in $\RR^n$.

Let us consider a change of the isothermal coordinates which in
complex form is given by: $t=t(s)$. Since the isothermal
coordinates are preserved, then the transformation $t=t(s)$ is
either holomorphic or antiholomorphic. Denote by $\tilde\Phi(s)$
the function with respect to the new coordinates $s$.

\textbf{The holomorphic case}. Using the definition
\eqref{Phi_def} we have:
\begin{equation*}
\tilde\Phi(s)=2\frac{\partial\x}{\partial s}=
2\frac{\partial\x}{\partial t}\frac{\partial t}{\partial s} \; .
\end{equation*}
This means that under a holomorphic change of the coordinates
$t=t(s)$ we have:
\begin{equation}\label{Phi_s-hol}
\tilde\Phi(s)=\Phi(t(s)) \frac{\partial t}{\partial s} \; .
\end{equation}

\textbf{The antiholomorphic case:} As in the above we find:
\begin{equation*}
\tilde\Phi(s)=2\frac{\partial\x}{\partial s}=
2\frac{\partial\x}{\partial \bar t}\frac{\partial \bar t}{\partial
s} \; .
\end{equation*}
Therefore under an antiholomorphic change of the coordinates
$t=t(s)$ we have:
\begin{equation}\label{Phi_s-antihol}
\tilde\Phi(s)=\bar\Phi(t(s)) \frac{\partial \bar t}{\partial s} \;
.
\end{equation}
In particular, under the change $t=\bar s$, the function $\Phi$ is
transformed in the following way:
\begin{equation}\label{Phi_s-t_bs}
\tilde\Phi(s)=\bar\Phi(\bar s) \; .
\end{equation}

Now, let us consider two surfaces $(\M,\x)$ and $(\hat\M,\hat\x)$
in $\RR^n$, parameterized by isothermal coordinates $t=u+\ii v$
defined in one and the same domain $\D \subset \CC$. Suppose that
$(\hat\M,\hat\x)$ is obtained from $(\M,\x)$ via a motion in
$\RR^n$ by the formula:
\begin{equation}\label{hat_M-M-mov}
\hat\x(t)=A\x(t)+b; \qquad A \in \mathbf{O}(n,\RR), \ b \in \RR^n
\; .
\end{equation}
Differentiating \eqref{hat_M-M-mov} we get the relation between
$\Phi$ and $\hat\Phi$:
\begin{equation}\label{hat_Phi-Phi-mov}
\hat\Phi(t)=A\Phi(t); \qquad A \in \mathbf{O}(n,\RR) \; .
\end{equation}
Conversely, if $\Phi$ and $\hat\Phi$ are related by
\eqref{hat_Phi-Phi-mov}, then we have $\hat\x_u=A\x_u$ and
$\hat\x_v=A\x_v$ which imply  \eqref{hat_M-M-mov} and
\eqref{hat_Phi-Phi-mov} are equivalent.

\section{Characterizing of minimal surfaces in $\RR^n$ by $\Phi$ }

Let $\M$ be a surface in $\RR^n$, parameterized by isothermal
coordinates and $\Phi$ is the function given by \eqref{Phi_def}.

Differentiating \eqref{Phi2}, we find:
\begin{equation}\label{Phi.dPhi_dbt}
\Phi\cdot\frac{\partial\Phi}{\partial\bar t}=0 \; .
\end{equation}
Since $\ds\frac{\partial\Phi}{\partial\bar t}$ is real, then it
follows that
\begin{equation}\label{bPhi.dPhi_dbt}
\bar\Phi\cdot\frac{\partial\Phi}{\partial\bar t}=0 \; .
\end{equation}
Equalities \eqref{Phi.dPhi_dbt} and \eqref{bPhi.dPhi_dbt} imply
that $\ds\frac{\partial\Phi}{\partial\bar t}$ is orthogonal to
$T(\M)$ and $\ds{\frac{\partial\Phi}{\partial\bar t}} \in N(\M)$.

Taking into account the last property and \eqref{dPhi_dbt} we
calculate:
\begin{equation*}
\begin{array}{rl}\ds
\frac{\partial\Phi}{\partial\bar t}\!\! &=
\ds\left(\frac{\partial\Phi}{\partial\bar t}\right)^\bot=
\frac{1}{2}(\Delta \x )^\bot= \frac{1}{2}(\x_{uu}+\x_{vv} )^\bot=
\frac{1}{2}(\nabla_{\x_u} \x_u + \nabla_{\x_v} \x_v )^\bot\\[2.5ex]
&=\ds\frac{1}{2}(\sigma(\x_u,\x_u) + \sigma(\x_v,\x_v) )= E\;
\frac{1}{2}(\sigma(\vX_1,\vX_1) + \sigma(\vX_2,\vX_2) ) = E\vH\; .
\end{array}
\end{equation*}
Thus we have:
\begin{equation}\label{dPhi_dbt-Delta_x-EH}
\frac{\partial\Phi}{\partial\bar t}=\frac{1}{2}\Delta \x = E\vH\;
.
\end{equation}

These equalities imply the following statement.
\begin{prop}
Let $(\M: \; (u,v) \to \x(u, v); \; (u,v) \in \D)$ be a surface in
$\RR^n$ parameterized by isothermal coordinates and $\Phi(t)$ be
the complex function:
\[\Phi(t)=2\frac{\partial\x}{\partial t}=\x_u-\ii\x_v ; \quad t=u+\ii v.\]

The following conditions are equivalent:
\begin{enumerate}
    \item the function $\Phi(t)$ is holomorphic $\left( \ds\frac{\partial\Phi}{\partial\bar t}= 0 \right)$;
    \vskip 1mm
    \item the function $\x (u,v)$ is harmonic $(\Delta \x = 0)$;
\vskip 2mm
    \item $(\M,\x)$ is a minimal surface in $\RR^n$ $(\vH=0)$.
\end{enumerate}
\end{prop}
\vskip 2mm Let $(\M,\x)$ be a minimal surface. We can introduce
the harmonic conjugate function $\y$ to $\x$ determined by the
conditions:
\[ \y_u=-\x_v; \quad \y_v=\x_u \ .\]
Then the function
\[\Psi=\x+\ii\y ,\]
is holomorphic and
\[\x=\Re\Psi ; \qquad \Phi=\x_u-\ii\x_v=\x_u+\ii\y_u=\frac{\partial\Psi}{\partial u}=\Psi' .\]
\smallskip
Since $\vH = 0$, we have:
\begin{equation}\label{sigma_22}
\sigma(\vX_2,\vX_2)=-\sigma(\vX_1,\vX_1).
\end{equation}
and
$$
\sigma(\x_v,\x_v)=E\sigma(\vX_2,\vX_2)=-E\sigma(\vX_1,\vX_1)=-\sigma(\x_u,\x_u).
$$
Then the formulas \eqref{dPhi_dt} and \eqref{dPhi_dt_bot} for the
derivative $\Phi'$  of $\Phi$ and its orthogonal projection on
$N_{p,C}(\M)$ become
\begin{equation}\label{PhiPr}
\Phi^\prime=\frac{\partial\Phi}{\partial
u}=\x_{uu}-\ii\x_{uv};\quad \Phi^{\prime \bot}=\x_{uu}^\bot -\ii
x_{uv}^\bot =\sigma (\x_u,\x_u)-\ii\sigma (\x_u,\x_v).
\end{equation}

Taking into account \eqref{sigma_uu_vv_uv}, we express $\sigma
(\x_u,\x_u)$, $\sigma (\x_v,\x_v)$ and $\sigma (\x_u,\x_v)$ by
means of $\Phi$:
\begin{equation}\label{sigma_uu_uv}
\begin{array}{l}
\sigma(\x_u,\x_u)=\ \ \,\Re (\Phi^{\prime \bot})=\ \
\,\ds\frac{1}{2}(\Phi^{\prime \bot}+\overline{\Phi^{\prime
\bot}})=
\ \ \,\ds\frac{1}{2}(\Phi^{\prime \bot}+{\overline{\Phi^\prime}}^\bot)\\[4mm]

\sigma(\x_v,\x_v)=     -\Re (\Phi^{\prime
\bot})=-\ds\frac{1}{2}(\Phi^{\prime \bot}+\overline{\Phi^{\prime
\bot}})=
-\ds\frac{1}{2}(\Phi^{\prime \bot}+{\overline{\Phi^\prime}}^\bot)\\[4mm]

\sigma(\x_u,\x_v)=-\Im (\Phi^{\prime \bot})=
\ds\frac{-1}{2\ii}(\Phi^{\prime \bot}-\overline{\Phi^{\prime
\bot}})= \ \ \!\,\ds\frac{\ii}{2}(\Phi^{\prime
\bot}-{\overline{\Phi^\prime}}^\bot).
\end{array}
\end{equation}

\section{Formulas for the Gauss curvature and the normal curvature}

Let $\M: \; (u,v) \to \x(u, v); \; (u,v) \in \D$ be a minimal
surface in $\RR^4$ parameterized by isothermal coordinates.
Suppose that $\n_1, \, \n_2$ be an orthonormal pair of normal
vector fields of $\M$, so that the quadruple
$\{\vX_1,\vX_2,\n_1,\n_2\}$ is right oriented in $\RR^4$. For any
normal vector $\n$ we denote by $A_{\n}$ the Weingarten operator
in $T(\M)$. This operator is connected with the second fundamental
form $\sigma$ by means of the equality: $A_{\n}\vX\cdot \vY=\sigma
(\vX,\vY)\cdot \n $. The condition $\vH=0$ implies that $\trace
A_{\n}=0 $ for any normal $\n$. Then the matrix representation of
the operators $A_{\n_1}$ and $A_{\n_2}$ has the following form:
\begin{equation}\label{A1A2_nu_lambda_rho_mu}
A_{\n_1}= \left(
\begin{array}{rr}
\nu      &  \lambda\\
\lambda  & -\nu
\end{array}
\right); \qquad A_{\n_2}= \left(
\begin{array}{rr}
\rho &  \mu\\
\mu  & -\rho
\end{array}
\right)
\end{equation}
Therefore
\begin{equation}\label{sigma_nu_lambda_rho_mu}
\begin{array}{l}
\sigma (\vX_1,\vX_1)=(\sigma (\vX_1,\vX_1)\cdot \n_1)\n_1+(\sigma
(\vX_1,\vX_1)\cdot \n_2)\n_2=
\nu \n_1 + \rho \n_2,\\
[1mm] \sigma (\vX_1,\vX_2)=(\sigma (\vX_1,\vX_2)\cdot
\n_1)\n_1+(\sigma (\vX_1,\vX_2)\cdot \n_2)\n_2=
\lambda \n_1 + \mu \n_2,\\
[1mm] \sigma (\vX_2,\vX_2)=-\sigma (\vX_1,\vX_1)=-\nu \n_1 - \rho
\n_2
\end{array}
\end{equation}

Denoting by $R$ the curvature tensor of the surface ${\M}$, the
Gauss equation and \eqref{sigma_22} imply that the Gauss curvature
$K$ of $\M$ is given by
\begin{equation}\label{K}
\begin{array}{rl}
K &=R(\vX_1,\vX_2,\vX_1,\vX_2)=R(\vX_1,\vX_2)\vX_2\cdot \vX_1\\
    &=-\sigma^2(\vX_1,\vX_1)-\sigma^2 (\vX_1,\vX_2)\;.
\end{array}
\end{equation}
On the other hand \eqref{sigma_nu_lambda_rho_mu} and \eqref{K}
imply the formula
\begin{equation}\label{K_nu_lambda_rho_mu}
K=-(\nu^2+\rho^2)-(\lambda^2+\mu^2)=-\nu^2-\lambda^2-\rho^2-\mu^2=\det(A_{\n_1})+\det(A_{\n_2})\;.
\end{equation}
In view of \eqref{PhiPr} we find the relation:
\begin{equation}\label{PhiPr_X}
\Phi^{\prime \bot}= E(\sigma (\vX_1,\vX_1)-\ii\sigma
(\vX_1,\vX_2)) \; .
\end{equation}
Thus we have:
$$
\begin{array}{rl}{\|\Phi^{\prime \bot}\|}^2=\Phi^{\prime \bot}\cdot\overline{\Phi^{\prime \bot}}
&=E^2(\sigma^2(\vX_1,\vX_1)+\sigma^2 (\vX_1,\vX_2)).\end{array}
$$
The last formula and \eqref{EG} give that
\begin{equation}\label{s2+s2}
\sigma^2(\vX_1,\vX_1)+\sigma^2 (\vX_1,\vX_2)=\frac{{\|\Phi^{\prime
\bot}\|}^2}{E^2}=\frac{4{\|\Phi^{\prime \bot}\|}^2}{\|\Phi\|^4}\;.
\end{equation}
Now \eqref{K} and \eqref{s2+s2} imply that
\begin{equation}\label{K_Phi}
K= \ds\frac{-4{\|\Phi^{\prime \bot}\|}^2}{\|\Phi\|^4}.
\end{equation}

We shall give to \eqref{K_Phi} another useful form. First we note
that the vector functions $\Phi$ and $\bar\Phi$ are orthogonal
with respect to the Hermitian dot product in ${\CC}^4$ and form an
orthogonal tangential basis. Therefore the tangential component of
$\Phi^\prime$ is given by

$$\Phi^{\prime\top}
=\ds\frac{\Phi^{\prime\top}\cdot\bar
\Phi}{\|\Phi\|^2}\Phi+\ds\frac{\Phi^{\prime\top}\cdot \Phi}{\|\bar
\Phi\|^2}\bar \Phi =\ds\frac{\Phi' \cdot \bar
\Phi}{\|\Phi\|^2}\Phi + \ds\frac{\Phi' \cdot \Phi}{\|\bar
\Phi\|^2}\bar \Phi. $$

Differentiating $\Phi^2=0$, we find $\Phi\cdot\Phi^\prime=0$. Then
we obtain for the projections of $\Phi'$ the following expression:
\begin{equation}\label{Phipn}
\Phi^{\prime\top}= \ds\frac{\Phi' \cdot \bar
\Phi}{\|\Phi\|^2}\Phi; \quad\quad
\Phi^{\prime\bot}=\Phi'-\Phi^{\prime\top}= \Phi'-\ds\frac{\Phi'
\cdot \bar \Phi}{\|\Phi\|^2}\Phi.
\end{equation}
Using a complex conjugation in \eqref{Phipn} we get:
\begin{equation*}
\begin{array}{rl}
{\|\Phi^{\prime\bot}\|}^2&=\Phi^{\prime\bot}\cdot\overline{\Phi^{\prime\bot}}=
\ds\frac{\|\Phi\|^2\|\Phi'\|^2-|\bar \Phi \cdot
\Phi'|^2}{\|\Phi\|^2}\ .
\end{array}
\end{equation*}
Since the bi-vector $\Phi\wedge\Phi'$ satisfies the equality
\[\|\Phi\wedge\Phi'\|^2=\|\Phi\|^2\|\Phi'\|^2-|\bar \Phi \cdot \Phi'|^2\]
then we have:
\begin{equation*}
\|\Phi^{\prime\bot}\|^2=\ds\frac{\|\Phi\|^2\|\Phi'\|^2-|\bar \Phi
\cdot \Phi'|^2}{\|\Phi\|^2}=
\ds\frac{\|\Phi\wedge\Phi'\|^2}{\|\Phi\|^2}\ .
\end{equation*}
Replacing into \eqref{K_Phi} we obtain:

\begin{equation}\label{K_Phi_bv}
K= \ds\frac{-4{\|\Phi^{\prime \bot}\|}^2}{\|\Phi\|^4}=
   \ds\frac{-4\|\Phi\wedge\Phi'\|^2}{\|\Phi\|^6}\ .
\end{equation}
\medskip
Further we find a similar formula for the normal curvature
$\varkappa$ of $\M$.

Denoting by $R^N$ the curvature tensor of the normal connection on
$\M$ we have:
\begin{equation}\label{kappa}
\begin{array}{rl}
\varkappa &= R^N(\vX_1,\vX_2,\n_1,\n_2)=A_{\n_1}\vX_1\cdot A_{\n_2}\vX_2-A_{\n_2}\vX_1\cdot A_{\n_1}\vX_2\\
&=2\nu\mu - 2\rho\lambda \; .
\end{array}
\end{equation}

Let us denote by $\det (\va,\vb,\vc,\vd)$ the determinant formed
by the coordinates of the four vectors $\va$, $\vb$, $\vc$ and
$\vd$, with respect to the standard basis in ${\CC}^4$. Using
\eqref{sigma_nu_lambda_rho_mu} we get
$$
\begin{array}{rl}
\det (\x_u,\x_v,\sigma(\x_u,\x_u),\sigma(\x_u,\x_v))&=E^3 \det (\vX_1,\vX_2,\sigma(\vX_1,\vX_1),\sigma(\vX_1,\vX_2))\\
&=E^3 \det (\vX_1,\vX_2,\nu\n_1,\mu\n_2)+E^3 \det (\vX_1,\vX_2,\rho\n_2,\lambda\n_1)\\
&=E^3(\nu\mu - \rho\lambda)\det (\vX_1,\vX_2,\n_1,\n_2)=E^3(\nu\mu
- \rho\lambda)$$
\end{array}
$$
Hence
\begin{equation}\label{numu2}
\nu\mu - \rho\lambda = \ds \frac{1}{E^3}\det
(\x_u,\x_v,\sigma(\x_u,\x_u),\sigma(\x_u,\x_v)).
\end{equation}

In the last equality we replace $\x_u$ and $\x_v$ taking into
account \eqref{xuxv} and find:
\begin{equation}\label{det1}
\begin{array}{l}
\ds \det (\x_u,\x_v,\sigma(\x_u,\x_u),\sigma(\x_u,\x_v))=\frac{\ii}{4}
\det (\Phi+\bar\Phi,\Phi-\bar\Phi,\sigma(\x_u,\x_u),\sigma(\x_u,\x_v))\\
\ds =-\frac{\ii}{2}\det
(\Phi,\bar\Phi,\sigma(\x_u,\x_u),\sigma(\x_u,\x_v)).
\end{array}
\end{equation}
In view of \eqref{sigma_uu_uv}, replacing $\sigma(\x_u,\x_u)$ and
$\sigma(\x_u,\x_v)$ we have:
\begin{equation}\label{det2}
\begin{array}{l}
\ds \det
(\Phi,\bar\Phi,\sigma(\x_u,\x_u),\sigma(\x_u,\x_v))=-\frac{\ii}{2}\det
(\Phi,\bar\Phi,\Phi^{\prime \bot},{\overline{\Phi^\prime}}^\bot).
\end{array}
\end{equation}
Now \eqref{det2} and \eqref{det1} imply that:
\begin{equation}\label{det3}
\begin{array}{l}
\ds \det (\x_u,\x_v,\sigma(\x_u,\x_u),\sigma(\x_u,\x_v))
=-\frac{1}{4}\det (\Phi,\bar\Phi,\Phi^{\prime
\bot},{\overline{\Phi^\prime}}^\bot)= -\frac{1}{4}\det
(\Phi,\bar\Phi,\Phi^\prime,\overline{\Phi^\prime}).
\end{array}
\end{equation}
Now \eqref{kappa},  \eqref{numu2} and \eqref{det3} give:
\begin{equation*}
\begin{array}{l}
\ds \varkappa = 2\nu\mu - 2\rho\lambda = \frac{2}{E^3}\det
(\x_u,\x_v,\sigma(\x_u,\x_u),\sigma(\x_u,\x_v))=
-\frac{1}{2E^3}\det
(\Phi,\bar\Phi,\Phi^\prime,\overline{\Phi^\prime}).
\end{array}
\end{equation*}
Finally, in view of \eqref{EG} we obtain the following formula for
$\varkappa$\::
\begin{equation}\label{kappa2}
\begin{array}{l}
\ds \varkappa = -\frac{4}{\|\Phi\|^6}\det
(\Phi,\bar\Phi,\Phi^\prime,\overline{\Phi^\prime}).
\end{array}
\end{equation}
Thus we obtained the following statement:
\begin{thm} The Gauss curvature $K$ and the normal curvature $\varkappa$ of any minimal surface
$(\M,\x)$ in $\RR^4$ parameterized by isothermal coordinates, are
given by the following formulas:

\begin{equation}\label{K_kappa_Phi}
K= \ds\frac{-4{\|\Phi^{\prime \bot}\|}^2}{\|\Phi\|^4}
=\ds\frac{-4\|\Phi\wedge\Phi'\|^2}{\|\Phi\|^6}; \quad\quad
\varkappa = -\ds\frac{4}{\|\Phi\|^6}\det
(\Phi,\bar\Phi,\Phi^\prime,\overline{\Phi^\prime}).
\end{equation}
\end{thm}

\section{Canonical coordinates on minimal surfaces in $\RR^4$. }

Let $\M$ be a surface in $\RR^4$. A point $p\in \M$ is said to be
super-conformal if the ellipse of curvature of $\M$ at the point
$p$ is a circle.

Now let $(\M ,\ \x=\Re\Psi)$ be a minimal surface in $\RR^4$
parameterized by isothermal coordinates $(u,v)$. A point $p\in \M$
is superconformal if
\begin{equation}\label{sconf_sigma}
\begin{array}{l}
\sigma (\vX_1,\vX_1)\bot \: \sigma (\vX_1,\vX_2)\\[2mm]
\sigma^2(\vX_1,\vX_1)=\sigma^2 (\vX_1,\vX_2)
\end{array}
\end{equation}
Next we express the condition \eqref{sconf_sigma} by means of the
function $\Phi$. Taking the square in \eqref{PhiPr_X}, we find:
\begin{equation}\label{PhiPr_bot2}
\begin{array}{rl}
{\Phi^{\prime \bot}}^2 \!\!\! &=E^2(\sigma^2(\vX_1,\vX_1)-\sigma^2
(\vX_1,\vX_2))-\ii\, 2E^2 \sigma (\vX_1,\vX_1)\sigma
(\vX_1,\vX_2).
\end{array}
\end{equation}
Comparing \eqref{sconf_sigma} with \eqref{PhiPr_bot2} we get the
equivalence
\begin{equation}\label{sconf}
\begin{array}{l}
\sigma (\vX_1,\vX_1)\bot \: \sigma (\vX_1,\vX_2)\\[2mm]
\sigma^2(\vX_1,\vX_1)=\sigma^2 (\vX_1,\vX_2)
\end{array}
\quad \Leftrightarrow \quad {\Phi^{\prime \bot}}^2=0
\end{equation}
Squaring the second equality of \eqref{Phipn}, we find:
$$
{\Phi^{\prime\bot}}^2 = {\Phi'}^2-2\Phi'\ds\frac{\Phi' \cdot \bar
\Phi}{\|\Phi\|^2}\Phi+\left(\ds\frac{\Phi' \cdot \bar
\Phi}{\|\Phi\|^2}\right)^2 \Phi^2.
$$
Taking into account $\Phi^2=0$ and $\Phi\cdot\Phi^\prime=0$, we
obtain:
\begin{equation}\label{PhiPr_bot2=PhiPr2}
{\Phi^{\prime \bot}}^2={\Phi^\prime}^2 \; .
\end{equation}
Thus we obtained the following proposition.
\begin{prop}
A point $p \in \M$ is superconformal if and only if $\Phi^{\prime
2}=0$.
\end{prop}
Now the fact that ${\Phi^\prime}^2$ is holomorphic implies the
following assertion.

\begin{thm}
If $\M$ is a connected minimal surface in $\RR^4$, then the set of
the superconformal points of $\M$ is either $\M$ or mostly a
countable set without limit points.
\end{thm}

Further we only consider minimal surfaces in $\RR^4$ without
superconformal points and call them \emph{minimal surfaces of
general type}. Any minimal surface of general type admits special
isothermal coordinates \cite{Itoh, TG, GM1}, such that the
coordinate vectors $\sigma (\vX_1,\vX_1)$ and $\sigma
(\vX_1,\vX_2)$ are directed along the principal axes of the
ellipse of curvature at the corresponding point. This means that
$\sigma (\vX_1,\vX_1)\bot \: \sigma (\vX_1,\vX_2)$. These
coordinates become uniquely determined adding the normalizing
condition $E^2(\sigma^2(\vX_1,\vX_1)-\sigma^2 (\vX_1,\vX_2))=\pm
1$. The sign "$+$" in the last formula corresponds to the case
when $\sigma^2(\vX_1,\vX_1)$ is directed along the major axis,
while the sign "$-$" corresponds to the case when
$\sigma^2(\vX_1,\vX_1)$ is directed along the minor axis of the
ellipse. We call the so described special isothermal coordinates
briefly \emph{canonical coordinates of the first type} and
\emph{canonical coordinates of the second type}.

In view of \eqref{PhiPr_bot2} we conclude that the isothermal
coordinates $(u,v)$ are canonical of the first kind if and only if
\begin{equation}\label{can1}
\begin{array}{l}
\sigma (\vX_1,\vX_1)\bot \: \sigma (\vX_1,\vX_2)\\[2mm]
E^2(\sigma^2(\vX_1,\vX_1)-\sigma^2 (\vX_1,\vX_2))=1
\end{array}
\quad \Leftrightarrow \quad {\Phi^\prime}^2={\Phi^{\prime
\bot}}^2=1
\end{equation}
The isothermal coordinates $(u,v)$ are canonical of the second
type if and only if

\begin{equation}\label{can2}
\begin{array}{l}
\sigma (\vX_1,\vX_1)\bot \: \sigma (\vX_1,\vX_2)\\[2mm]
E^2(\sigma^2(\vX_1,\vX_2)-\sigma^2 (\vX_1,\vX_1))=1
\end{array}
\quad \Leftrightarrow \quad {\Phi^\prime}^2={\Phi^{\prime
\bot}}^2=-1
\end{equation}

Using the properties of the function $\Phi$, we shall show that
any minimal surface of general type in $\RR^4$ carries locally
canonical coordinates of both types.

Let $(u,v)$ be isothermal coordinates on $\M$ and denote $t=u+vi$.
Consider the change $t=t(\tilde t\:)$, where $\tilde t\:$ is a new
complex coordinate. Next we find the conditions under which the
change $t(\tilde t\,)$ gives canonical coordinates. Firstly, the
new coordinates $\tilde t\:$ have to be isothermal. Therefore the
transformation $t=t(\tilde t\,)$ is conformal in ${\CC}$, which
means that $t(\tilde t\,)$ is either a holomorphic or an
antiholomorphic function. It is enough to consider only the case
of a holomorphic change $t=t(\tilde t\,)$.

Let $\tilde\Psi$ be the holomorphic function representing ${\M}$
with respect to the new coordinates, and $\tilde\Phi$ be its
derivative. Then we have:
\begin{equation}\label{tildPhi}
\tilde\Phi = \tilde\Psi'_{\tilde t} = \Psi'_t t' = \Phi t'
\end{equation}
Further we find: $\tilde\Phi'_{\tilde t}=\Phi'_t t'^2+\Phi t''$.
Since $\Phi$ is tangent to ${\M}$, then $\Phi^\bot=0$ and
therefore:

\begin{equation}\label{tildPhiPr2}
\begin{array}{lll}
\tilde\Phi_{\tilde t}'^\bot &=& (\Phi'_t t^{\prime \, 2}+\Phi t'')^\bot = \Phi_t'^\bot t^{\prime \, 2};\\
\left.\tilde\Phi_{\tilde t}^{\prime\bot}\right.^2 &=&
{\Phi_t^{\prime \bot}}^2 t'^4.
\end{array}
\end{equation}

According to \eqref{can1} and \eqref{can2} the change $\tilde t$
determines canonical coordinates if $\left.\tilde\Phi_{\tilde
t}^{\prime\bot}\right.^2=\pm 1$. Equalities \eqref{tildPhiPr2}
imply that if ${\Phi_t^{\prime \bot}}^2=0$, then
$\left.\tilde\Phi_{\tilde t}^{\prime\bot}\right.^2=0$, which is
the condition $\M$ to be superconformal. Hence, there do not exist
canonical coordinates on a superconformal surface.

If ${\M}$ is a minimal surface of general type, i.e.
${\Phi^{\prime \bot}}^2 \neq 0$, then $\tilde t$ determines
canonical coordinates if ${\Phi_t^{\prime \bot}}^2 t'^4=\pm 1$.
Thus the function $t(\tilde t\,)$ satisfies the following first
order ordinary differential equation:
\begin{equation}\label{eqcan}
\sqrt[4]{\pm{\Phi_t^{\prime \bot}}^2}\:dt = d{\tilde t}
\end{equation}
According to \eqref{PhiPr_bot2=PhiPr2} the left hand side of
\eqref{eqcan} is holomorphic and  after integrating of
\eqref{eqcan} we obtain $\tilde t$ as a holomorphic function of
$t$.

The condition ${\Phi_t^{\prime \bot}}^2 \neq 0$ means that
${\tilde t}'\neq 0$ and therefore the correspondence between
${\tilde t}$ and $t$ is one to one. Hence ${\tilde t}$ determines
new isothermal coordinates satisfying $\left.\tilde\Phi_{\tilde
t}^{\prime\bot}\right.^2=\pm 1$, i.e. the new coordinates are
canonical.

Thus we proved the following assertion.
\begin{prop}
Any minimal surface ${\M}$ in $\RR^4$ of general type admits
locally canonical coordinates of the first or of the second type.
\end{prop}

Next we consider the question of uniqueness of the canonical
coordinates.

 Let us assume that $t$ and $\tilde t$ are canonical coordinates on ${\M}$ of one and the same type.
Then $t = t(\tilde t\,)$ is either holomorphic or antiholomorphic
function.

First we consider the holomorphic case. Then the conditions
\eqref{can1}, \eqref{can2} and \eqref{tildPhiPr2} imply the
equalities:
$$
\pm 1 = \left.\tilde\Phi_{\tilde
t}^{\prime\bot}\right.^2={\Phi_t^{\prime \bot}}^2 t'^4 = \pm 1
t'^4 = \pm t'^4
$$
Therefore $t'^4 = 1$ and hence $t' = \pm 1;\ \pm i$. This implies
that $t$ and $\tilde t$ satisfy one of the following relations:
$t=\pm\tilde t+c;\ \pm i\tilde t+c$, where $c=\text{const}$.

The antiholomorphic case reduces to the previous case by the
change ${\tilde t}=\bar s$. From the last equality it follows that
$t=\pm\bar{\tilde t}+c;\ \pm i\bar{\tilde t}+c$. The last eight
relations mean that the canonical coordinates of one and the same
type are unique up to numbering and change of the direction of the
coordinate lines.

Finally, let us consider the relation between the canonical
coordinates of different type. Let $t=u+vi$ be canonical
coordinates of the first type an let us introduce new coordinates
by means of $t = e^{\frac{\pi i}{4}}\tilde t$. We find from here
that $t'^4=-1$. Taking into account \eqref{tildPhiPr2} we obtain
that $\left.\tilde\Phi_{\tilde t}^{\prime\bot}\right.^2=-1$ and
hence $\tilde t$ determines canonical coordinates of the second
type. Geometrically this means that the canonical coordinates of
both types are related to each other by a rotation to an angle
$\frac{\pi}{4}$ in the coordinate plane $(u,v)$

Let $(\M,\x)$ be a minimal surface of general type in $\RR^4$
parameterized by canonical coordina\-tes of the first type. Up to
now the vectors $\n_1$ and $\n_2$ were only an orthonormal pair in
$N(\M)$. If the coordinates are canonical, then $\sigma
(\vX_1,\vX_1)\bot \: \sigma (\vX_1,\vX_2)$. Therefore we can
choose $\n_1$ and $\n_2$ to have the directions of $\sigma
(\vX_1,\vX_1)$ and $\sigma (\vX_1,\vX_2)$, i.e. along the
principal axes of the ellipse of curvature at the corresponding
point. More precisely, let $\n_1$ be the unit normal vector with
the direction of $\sigma (\vX_1,\vX_1)$, and $\n_2$ be the unit
normal vector such that the quadruple $(\vX_1,\vX_2,\n_1,\n_2)$
determine a positive oriented orthonormal basis in $\RR^4$. Then
$\n_2$ is collinear with $\sigma (\vX_1,\vX_2)$. Under these
conditions formulas \eqref{sigma_nu_lambda_rho_mu} become
\begin{equation}\label{sigma_nu_mu}
\begin{array}{l}
\sigma (\vX_1,\vX_1)=\phantom{-} \nu \n_1 \\
\sigma (\vX_1,\vX_2)=\phantom{-} \mu \n_2 \\
\sigma (\vX_2,\vX_2)=         -  \nu \n_1
\end{array}\; ;
\qquad \nu>0 \; .
\end{equation}
which means that $\lambda=0$, $\rho=0$ and formulas
\eqref{A1A2_nu_lambda_rho_mu} become as follows:
\begin{equation}\label{A1A2_nu_mu}
A_{\n_1}= \left(
\begin{array}{rr}
\nu  &  0\\
0    & -\nu
\end{array}
\right), \qquad A_{\n_2}= \left(
\begin{array}{rr}
0    & \mu\\
\mu  & 0
\end{array}
\right).
\end{equation}
The functions $\nu$ and $\mu$ satisfy the following relations:
\begin{equation}\label{numu1}
\begin{array}{ll}
\nu \!\! & = \|\sigma (\vX_1,\vX_1)\|\\
|\mu | \!\! & = \|\sigma (\vX_1,\vX_2)\|
\end{array}\; ;
\qquad \nu > |\mu |\; .
\end{equation}
These functions do not depend on the canonical coordinates and are
invariants of a minimal surface in $\RR^4$ \cite{GM1}. According
to \eqref{sigma_nu_mu}, these functions determine completely the
second fundamental form of $\M$. The second condition in
\eqref{can1} implies that the first fundamental form is also
completely determined by the formula:

\begin{equation}\label{E_nu_mu}
E=G=\frac{1}{\sqrt{\nu^2-\mu^2}} \; .
\end{equation}

Next we obtain explicit formulas expressing the pair $(\nu ,\mu)$
by the pair $(K ,\varkappa)$ and vice versa. Under the condition
$\lambda=0$ and $\rho=0$ formulas \eqref{K_nu_lambda_rho_mu} have
the following form \cite{GM1}:

\begin{equation}\label{K_kappa_nu_mu}
K= -\nu^2-\mu^2 < 0\; ; \quad\quad \varkappa = 2\nu\mu; \qquad -K
> |\varkappa|.
\end{equation}
Therefore
\begin{equation}\label{nu_mu_K_kappa}
\nu = \frac{1}{2}(\sqrt{-K+\varkappa}+\sqrt{-K-\varkappa}),\quad
\mu = \frac{1}{2}(\sqrt{-K+\varkappa}-\sqrt{-K-\varkappa}).
\end{equation}

Further we give formulas for $\nu$, $\mu$ and $\varkappa$, with
respect to canonical coordinates of the first type.

Taking into account \eqref{can1} and \eqref{EG} we have:
\begin{equation}\label{s2-s2}
\sigma^2(\vX_1,\vX_1)-\sigma^2
(\vX_1,\vX_2)=\frac{1}{E^2}=\frac{4}{\|\Phi\|^4}\;.
\end{equation}

From here and \eqref{s2+s2} we get:
\begin{equation}\label{nu2+-mu2}
\begin{array}{ll}
\nu^2+\mu^2 \!\! & = \ds \frac{4{\|\Phi^{\prime \bot}\|}^2}{\|\Phi\|^4}\\[4mm]
\nu^2-\mu^2 \!\! & = \ds \frac{4}{\|\Phi\|^4}
\end{array}
\quad \Leftrightarrow \quad
\begin{array}{ll}
\nu^2 \!\! & = \ds \frac{2({\|\Phi^{\prime \bot}\|}^2+1)}{\|\Phi\|^4}\\[4mm]
\mu^2 \!\! & = \ds \frac{2({\|\Phi^{\prime
\bot}\|}^2-1)}{\|\Phi\|^4}.
\end{array}
\end{equation}

In view of \eqref{nu2+-mu2} we find
\begin{equation}\label{modkappa}
|\varkappa \:|=|2\nu \mu |= \ds \frac{4\sqrt{{\|\Phi^{\prime
\bot}\|}^4-1}}{\|\Phi\|^4}.
\end{equation}

\section{Weierstrass representations for minimal surfaces in $\RR^4$.}

First we give some Weierstrass representations for minimal
surfaces of general type parame\-terized by isothermal coordinates.
Such kind of formulas have been written by a number of
mathematicians: e.g. Eisenhart \cite{E}, Hoffman and Osserman \cite{HO}.

Let $(\M,\x)$:\; $\x=\Re\Psi$ be a minimal surfaces in $\RR^4$,
parameterized by isothermal coordinates and let
$\Phi=\Psi^\prime$. If $\Phi=(\phi_1,\phi_2,\phi_3,\phi_4)$, then
the condition for isothermal coordinates $\Phi^2=0$ has the form:

\begin{equation}\label{phi2coord}
\phi_1^2+\phi_2^2+\phi_3^2+\phi_4^2=0.
\end{equation}

This equality can be "parameterized"\, in different ways by means
of three holomorphic functions.

First we shall find a representation of $\Phi$ by means of
trigonometric functions. Equality \eqref{phi2coord} is equivalent
to one of the following equalities:
$$\phi_1^2+\phi_2^2=-\phi_3^2-\phi_4^2;\quad
  \phi_1^2+\phi_3^2=-\phi_2^2-\phi_4^2;\quad
  \phi_1^2+\phi_4^2=-\phi_2^2-\phi_3^2.$$

At least one of the functions $\phi_1^2+\phi_2^2$,
$\phi_1^2+\phi_3^2$ and $\phi_1^2+\phi_4^2$ has to be different
from zero. (The inverse leads by means of \eqref{phi2coord} to
$\phi_1^2=\phi_2^2=\phi_3^2=\phi_4^2=0$, which is impossible.)
Without loss of generality we can assume that
$\phi_1^2+\phi_2^2\neq 0$. Therefore, there exists a holomorphic
function $f \neq 0$, such that:
\begin{equation}\label{f}
f^2=\phi_1^2+\phi_2^2=-\phi_3^2-\phi_4^2.
\end{equation}

The last equality is equivalent to
\begin{equation}\label{f2}
\left(\frac{\phi_1}{f}\right)^2+\left(\frac{\phi_2}{f}\right)^2=
\left(\frac{\phi_3}{if}\right)^2+\left(\frac{\phi_4}{if}\right)^2=1.
\end{equation}
It follows from here that there exist holomorphic functions $h_1$
and $h_2$, such that
$$
\frac{\phi_1}{f}=\cos h_1;\quad \frac{\phi_2}{f}=\sin h_1; \quad
\frac{\phi_3}{if}=\cos h_2;\quad \frac{\phi_4}{if}=\sin h_2.
$$
Thus we obtain the following representation of the vector function
$\Phi$:
\begin{equation}\label{W1}
\Phi: \quad
\begin{array}{rlr}

\phi_1 &=& f\cos h_1,\\
[1mm]
\phi_2 &=& f\,\sin h_1,\\
[1mm]
\phi_3 &=& if\cos h_2,\\
[1mm]
\phi_4 &=& if\,\sin h_2.\\
\end{array}
\end{equation}
Hence, any minimal surface  $\rm M$ in $\RR^4$, parameterized by
isothermal parameters has a Weier\-strass representation of the
type \eqref{W1}.

Conversely, for any three holomorphic functions $(f \neq
0,h_1,h_2)$ determined in a region $\rm D \subset{\mathbb C}$,
formulas \eqref{W1} generate a holomorphic function $\Phi$ with
values in ${\mathbb C}^4$. The condition $f \neq 0$ gives $\Phi
\neq 0$. By direct calculations, formulas \eqref{W1} imply
\eqref{phi2coord}, which is $\Phi^2=0$. Determining $\Psi$ by the
condition $\Psi^\prime = \Phi$ and defining $\rm M:\; x={\rm
Re}(\Psi)$, we obtain a minimal surface $\rm M$ in $\RR^4$,
parameterized by isothermal coordinates.

Hence, any triplet of holomorphic functions $(f \neq 0,h_1,h_2)$
generates a minimal surface in $\RR^4$ via formulas \eqref{W1}.

Finally, we shall establish to what extent the triple $(f \neq
0,h_1,h_2)$ is determined by $\Phi$. For that purpose, let us
assume that one and the same function $\Phi$ is represented by
\eqref{W1} via two different triplets $(f \neq 0,h_1,h_2)$ and
$(\hat f \neq 0,\hat h_1,\hat h_2)$. It is seen from \eqref{f}
that, $f$ is determined by $\Phi$ up to a sign. Therefore, two
cases are possible. If $\hat f = f$, then $\hat h_1$ and $\hat
h_2$ differ from $h_1$ and $h_2$ by constants even multiples to
$\pi$. If $\hat f = -f$, then $\hat h_1$ and $\hat h_2$ differ
from $h_1$ and $h_2$ by constants odd multiples to $\pi$. Thus we
have:
\begin{equation*}
\begin{array}{ll}
\hat f   \!\! &=\: f\\
\hat h_1 \!\! &=\: h_1 + 2k_1\pi\\
\hat h_2 \!\! &=\: h_2 + 2k_2\pi
\end{array}\quad
\text{or}\quad
\begin{array}{ll}
\hat f   \!\! &=\: -f\\
\hat h_1 \!\! &=\: h_1 + (2k_1+1)\pi\\
\hat h_2 \!\! &=\: h_2 + (2k_2+1)\pi
\end{array};\quad \ \
\begin{array}{l}
k_1=\text{const}\\
k_2=\text{const}
\end{array}
\end{equation*}
\vspace {2mm}

Using \eqref{W1} we can obtain another forms of the Weierstrass
representation for minimal surfaces applying different
replacements.

In order to obtain the Weierstrass representation by means of
hyperbolic functions, we make the following replacements in
\eqref{W1}:
$$
f\rightarrow if; \quad h_1 \rightarrow -ih_1; \quad h_2
\rightarrow \pi +ih_2.
$$
Thus we obtain the following Weierstrass representation by means
of hyperbolic functions:
\begin{equation}\label{W2}
\Phi: \quad
\begin{array}{rlr}

\phi_1 &=& if\cosh h_1,\\
[1mm]
\phi_2 &=&  f\,\sinh h_1,\\
[1mm]
\phi_3 &=&  f\cosh h_2,\\
[1mm]
\phi_4 &=& if\,\sinh h_2.\\
\end{array}
\end{equation}

Let us introduce the functions $w_1$ и $w_2$ instead of $h_1$ и
$h_2$ in $\eqref{W2}$ as follows:
\begin{equation}\label{w1w2}
\begin{array}{l}
w_1=h_1+h_2,\\
w_2=h_1-h_2.
\end{array}
\end{equation}
Thus we obtain Weierstrass representation of the following type:
\begin{equation}\label{W5}
\Phi: \quad
\begin{array}{rlr}
\phi_1 &=& if \cosh \ds\frac{w_1+w_2}{2}\,,\\[4mm]
\phi_2 &=&  f \sinh \ds\frac{w_1+w_2}{2}\,,\\[4mm]
\phi_3 &=&  f \cosh \ds\frac{w_1-w_2}{2}\,,\\[4mm]
\phi_4 &=& if \sinh \ds\frac{w_1-w_2}{2}\,.\\
\end{array}
\end{equation}

Further, let us introduce the functions $g_1$ and $g_2$ by the
following formulas:
\begin{equation}\label{g}
g_1=e^{w_1}; \quad g_2=e^{w_2}.
\end{equation}
With the aid of these functions, in view of \eqref{W5}, we obtain
Weierstrass representation, which is a natural analogue of the
classical Weierstrass representation for minimal surfaces in
$\RR^3$.

First, we calculate $\phi_1$:
$$
\begin{array}{rl}
\phi_1 &= \ds\frac{if}{2}
(e^{\frac{w_1+w_2}{2}}+e^{-\frac{w_1+w_2}{2}})=
          \ds\frac{if}{2} e^{-\frac{w_1}{2}}e^{-\frac{w_2}{2}} (e^{w_1+w_2}+ 1)\\[6mm]
       &= \ds\frac{if}{2\sqrt{g_1 g_2}}(e^{w_1}e^{w_2}+ 1)=
          \ds\frac{if}{2\sqrt{g_1 g_2}}(g_1 g_2+1).
\end{array}
$$
Analogously to the above, we compute $\phi_2$:
$$
\phi_2 = \ds\frac{f}{2}
(e^{\frac{w_1+w_2}{2}}-e^{-\frac{w_1+w_2}{2}})=
         \ds\frac{f}{2\sqrt{g_1 g_2}}(g_1 g_2-1).
$$
In a similar way we find $\phi_3$:
$$
\begin{array}{rl}
\phi_3 &= \ds\frac{f}{2}
(e^{\frac{w_1-w_2}{2}}+e^{-\frac{w_1-w_2}{2}})=
          \ds\frac{f}{2} e^{-\frac{w_1}{2}}e^{-\frac{w_2}{2}} (e^{w_1}+e^{w_2})\\[6mm]
       &= \ds\frac{f}{2\sqrt{g_1 g_2}}(e^{w_1}+e^{w_2})=
          \ds\frac{f}{2\sqrt{g_1 g_2}}(g_1+g_2).
\end{array}
$$
Finally we calculate $\phi_4$:
$$
\phi_4 = \ds\frac{if}{2}
(e^{\frac{w_1-w_2}{2}}-e^{-\frac{w_1-w_2}{2}})=
         \ds\frac{if}{2\sqrt{g_1 g_2}}(g_1 - g_2).
$$
\noindent In the last four formulas, we make the change:
\begin{equation}\label{ff}
f \rightarrow f2\sqrt{g_1 g_2}
\end{equation}
and obtain the following polynomial Weierstrass representation:
\begin{equation}\label{W6}
\Phi: \quad
\begin{array}{rll}
\phi_1 &=& if(g_1 g_2+1),\\
\phi_2 &=& \ f(g_1 g_2-1),\\
\phi_3 &=& \ f(g_1+g_2),\\
\phi_4 &=& if(g_1-g_2).\\
\end{array}
\end{equation}

Conversely, if $(f \neq 0,g_1,g_2)$ are three holomorphic
functions, determined in a region in ${\mathbb C}$, then by virtue
of \eqref{W6} we obtain a holomorphic function $\Phi$ with values
in ${\mathbb C}^4$. It follows from $f \neq 0$ that $\Phi \neq 0$.
It is easy to see by direct calculations that \eqref{W6} implies
\eqref{phi2coord}, which is $\Phi^2=0$. Therefore, if we define
$\Psi$ by the equality $\Psi^\prime = \Phi$, then the surface
${\rm M}:\; x=\Re (\Psi)$, will be a minimal surface in $\RR^4$,
parameterized by isothermal coordinates. Hence any triplet of
holomorphic functions $(f \neq 0,g_1,g_2)$ generates a minimal
surface in $\RR^4$ via formulas \eqref{W6}.

Finally, we shall obtain that the triplet $(f \neq 0,g_1,g_2)$ is
determined uniquely by $\Phi$. For that purpose we express the
functions $f$, $g_1$ и $g_2$ explicitly by $\Phi$. As an immediate
consequence of \eqref{W6}, we find:
$$
\begin{array}{l}
i\phi_1+\phi_2=-f(g_1 g_2+1)+f(g_1 g_2-1)=-2f,\\
\phi_3+i\phi_4=f(g_1+g_2)-f(g_1-g_2)=2fg_2,\\
\phi_3-i\phi_4=f(g_1+g_2)+f(g_1-g_2)=2fg_1.
\end{array}
$$
The above equalities imply the following formulas for $f$, $g_1$ и
$g_2$:
\begin{equation}\label{fg1g2}
f=-\ds\frac{1}{2}(i\phi_1+\phi_2); \quad
g_1=-\ds\frac{\phi_3-i\phi_4}{i\phi_1+\phi_2}; \quad
g_2=-\ds\frac{\phi_3+i\phi_4}{i\phi_1+\phi_2}.
\end{equation}
\section{Canonical Weierstrass representations of minimal surfaces}

A Weierstrass representation with respect to isothermal
coordinates is said to be \emph{canonical of the first or the
second type} if the coordinates are in addition canonical of the
first or the second type, respectively. In this section we shall
only consider canonical Weierstrass representations of the first
type.

\subsection{Preliminary calculations}
In order to obtain canonical Weierstrass representations of
minimal surfaces in $\RR^4$ we give first some relations between
the functions $f$, $h_1$ and $h_2$, that are used in the
Weierstrass representation of minimal surfaces.

Here we prefer to use the representation \eqref{W2} via hyperbolic
functions. From now on we use the scalar holomorphic functions
$w_1$ and $w_2$, defined by \eqref{w1w2} and the vector
holomorphic function $a$ defined in the following way:
\begin{equation}\label{adef}
a=\ds\frac{\Phi}{f}
\end{equation}

Taking into account \eqref{W2} and \eqref{adef} we get the
following formulas for the functions $a$, $\bar a$, $a'$ and $\bar
{a'}$:

\begin{equation}\label{a}
\begin{array}{l}
a=(\ \ \, i\cosh h_1,\sinh h_1,\cosh h_2,\ \ i\sinh h_2),\\
\bar a=(-i\cosh \bar h_1,\sinh \bar h_1,\cosh \bar h_2,-i\sinh \bar h_2),\\
a'\! =(\ \ \, ih'_1 \sinh h_1,h'_1 \cosh h_1,h'_2 \sinh h_2,\ \ ih'_2 \cosh h_2),\\
\bar{a'}\! =(-i\bar{h'_1} \sinh \bar h_1,\bar{h'_1} \cosh \bar
h_1,\bar{h'_2} \sinh \bar h_2,-i\bar{h'_2} \cosh \bar h_2).
\end{array}
\end{equation}

Now we can find the inner products between $a$, $\bar a$, $a'$ and
$\bar {a'}$.

The condition $\Phi^2=0$, implies that $a^2=0$. By means of
differentiation and complex conjugation we get
\begin{equation}\label{a1}
a^2=aa'=\bar a^2=\bar a \bar {a'}=0
\end{equation}
Multiplying equations \eqref{a} we also find
\begin{equation}\label{a2}
\begin{array}{rl}
\|a\|^2=a\bar a &= \cosh h_1 \cosh \bar h_1 + \sinh h_1 \sinh \bar h_1 + \cosh h_2 \cosh \bar h_2 + \sinh h_2 \sinh \bar h_2\\
                                &= 2\cosh(\Re w_1)\cosh(\Re w_2);
\end{array}
\end{equation}
\begin{equation}\label{a3}
\begin{array}{rl}
a\bar {a'} &= \bar{h'_1}\cosh h_1 \sinh \bar h_1 + \bar{h'_1}\sinh
h_1 \cosh \bar h_1 +
              \bar{h'_2}\cosh h_2 \sinh \bar h_2 + \bar{h'_2}\sinh h_2 \cosh \bar h_2\\
                               &= \bar{h'_1}\sinh (2\Re h_1 ) + \bar{h'_2}\sinh (2\Re h_2 );
\end{array}
\end{equation}
\begin{equation}\label{a4}
 \bar a a' = \overline{a\bar {a'}} = h'_1 \sinh (2\Re h_1 ) + h'_2 \sinh (2\Re h_2 );
\end{equation}
\begin{equation}\label{a5}
\begin{array}{rl}
a'^2 &= -h'^2_1\sinh^2 h_1 + h'^2_1\cosh^2 h_1 + h'^2_2\sinh^2 h_2 - h'^2_2\cosh^2 h_2\\
     &= h'^2_1-h'^2_2 = w'_1 w'_2;
\end{array}
\end{equation}
\begin{equation}\label{a6}
\begin{array}{rl}
\|a'\|^2 = a'\bar {a'} &= |h'_1|^2 \sinh h_1 \sinh \bar h_1 + |h'_1|^2 \cosh h_1 \cosh \bar h_1\\
                       &+\  |h'_2|^2 \sinh h_2 \sinh \bar h_2 + |h'_2|^2 \cosh h_2 \cosh \bar h_2\\
                                                      &= |h'_1|^2\cosh (2\Re h_1 ) + |h'_2|^2\cosh (2\Re h_2 ).
\end{array}
\end{equation}

Next we find formulas for $a^{\prime\bot}$, ${a^{\prime\bot}}^2$
and ${\|a^{\prime\bot}\|}^2$ expressed by $h_1$, $h_2$ and by
$w_1$ и $w_2$, respectively. We have
$a^{\prime\bot}=a'-a^{\prime\top}$. The equality $a^2=0$ means
that the vectors $a$ and $\bar a$ are mutually orthogonal with
respect to the Hermitian inner product in ${\mathbb C}^4$.
Therefore $a^{\prime\top}$ being tangent to $\rm M$, can be
represented with respect to the orthogonal basis $(a,\bar a)$ in
the following way:
$$a^{\prime\top}=\ds\frac{a^{\prime\top}\cdot\bar a}{\|a\|^2}a+\ds\frac{a^{\prime\top}\cdot a}{\|\bar a\|^2}\bar a
                =\ds\frac{a' \cdot \bar a}{\|a\|^2}a + \ds\frac{a' \cdot a}{\|\bar a\|^2}\bar a. $$
In view of \eqref{a1} we get $a' \cdot a = 0$. Hence
\begin{equation}\label{apn}
a^{\prime\top}= \ds\frac{a' \cdot \bar a}{\|a\|^2}a; \quad\quad
a^{\prime\bot}=a'-a^{\prime\top}= a'-\ds\frac{a' \cdot \bar
a}{\|a\|^2}a.
\end{equation}
Squaring the second equality of \eqref{apn} we find
$$
{a^{\prime\bot}}^2 = {a'}^2-2a'\ds\frac{a' \cdot \bar
a}{\|a\|^2}a+\left(\ds\frac{a' \cdot \bar a}{\|a\|^2}\right)^2
a^2.
$$
According to \eqref{a1} $a' \cdot a = 0$ and $a^2=0$. Therefore we
have ${a^{\prime\bot}}^2 = {a'}^2$. By means of \eqref{a5} we
obtain the required expression for the function
${a^{\prime\bot}}^2$:
\begin{equation}\label{apn2}
{a^{\prime\bot}}^2 = {a'}^2 = {h'_1}^2-{h'_2}^2 = w'_1 w'_2
\end{equation}

By means of complex conjugation in \eqref{apn} we find the
following formula for ${\|a^{\prime\bot}\|}^2$:
\begin{equation}\label{mapn2}
\begin{array}{rl}
{\|a^{\prime\bot}\|}^2&=a^{\prime\bot}\cdot\overline{a^{\prime\bot}}=
\left(a'-\ds\frac{a' \cdot \bar a}{\|a\|^2}a\right)\left(\bar{a'}-\ds\frac{\bar{a'} \cdot a}{\|a\|^2}\bar a \right)\\[4mm]
                                                                  &=\|a'\|^2 - \ds\frac{|\bar{a'} \cdot a|^2}{\|a\|^2} - \ds\frac{|a' \cdot \bar a|^2}{\|a\|^2}+
                                              \ds\frac{|a' \cdot \bar a|^2}{\|a\|^4}\|a\|^2 = \|a'\|^2 - \ds\frac{|\bar{a'} \cdot a|^2}{\|a\|^2}\\[4mm]
                                            &=\ds\frac{\|a\|^2\|a'\|^2-|\bar a \cdot a'|^2}{\|a\|^2}
\end{array}
\end{equation}

Let us denote the numerator in formula \eqref{mapn2} by $k_1$.
Applying formulas \eqref{a2}, \eqref{a4} and \eqref{a6} after the
corresponding simplification we find:
\begin{equation}\label{k1}
\begin{array}{rl}
k_1 &= \|a\|^2\|a'\|^2-|\bar a \cdot a'|^2\\
    &= (|h'_1|^2+|h'_2|^2)(1+\cosh (2\Re h_1 )\cosh (2\Re h_2 ))\\
        &- 2\Re(h'_1\bar h'_2)\sinh (2\Re h_1 )\sinh (2\Re h_2 )
\end{array}
\end{equation}

Denoting the determinant of the vectors $a$, $\bar a$, $a'$ and
$\bar {a'}$ by $-k_2$, by direct calculations we find that
\begin{equation}\label{k2}
\begin{array}{rl}
k_2 &= -\det(a,\bar a , a' , \bar {a'})\\
    &= 2\Re(h'_1\bar h'_2)(1+\cosh (2\Re h_1 )\cosh (2\Re h_2 ))\\
        &- (|h'_1|^2+|h'_2|^2)\sinh (2\Re h_1 )\sinh (2\Re h_2 )
\end{array}
\end{equation}

Adding and subtracting equalities \eqref{k1} and \eqref{k2} we
obtain:
\begin{equation}\label{k1+k2}
k_1+k_2 = 2|h'_1+h'_2|^2 \cosh^2 (\Re h_1-\Re h_2 ).
\end{equation}
\begin{equation}\label{k1-k2}
k_1-k_2 = 2|h'_1-h'_2|^2 \cosh^2 (\Re h_1+\Re h_2 ).
\end{equation}

Equalities \eqref{k1+k2} and \eqref{k1-k2} give the following
expressions for $k_1$ and $k_2$:
\begin{equation}\label{k12h}
\begin{array}{l}
k_1 = |h'_1+h'_2|^2 \cosh^2 (\Re h_1-\Re h_2 ) + |h'_1-h'_2|^2 \cosh^2 (\Re h_1+\Re h_2 )\\
k_2 = |h'_1+h'_2|^2 \cosh^2 (\Re h_1-\Re h_2 ) - |h'_1-h'_2|^2
\cosh^2 (\Re h_1+\Re h_2 )
\end{array}
\end{equation}

Replacing $h_1$ and $h_2$ by $w_1$ and $w_2$, respectively, we
get:
\begin{equation}\label{k12w}
\begin{array}{l}
k_1 = |w'_1|^2 \cosh^2 (\Re w_2 ) + |w'_2|^2 \cosh^2 (\Re w_1 )\\
k_2 = |w'_1|^2 \cosh^2 (\Re w_2 ) - |w'_2|^2 \cosh^2 (\Re w_1 )
\end{array}
\end{equation}

\subsection{Canonical Weierstrass representations of minimal surfaces in $\RR^4$}

Let the minimal surface of general type ${\M}$ in $\RR^4$ be
parameterized by canonical coordinates of the first type and
assume that $\M$ is given by the representation \eqref{W2} by
means of hyperbolic functions. The condition \eqref{can1} for the
coordinates to be canonical implies a relation between the three
functions $f$, $h_1$ and $h_2$. In order to obtain this relation,
we express the condition ${\Phi^{\prime \bot}}^2=1$ via $f$, $h_1$
and $h_2$. In view of \eqref{adef} we have $\Phi=f\va$ and
therefore $\Phi'=f'\va+f\va'$. Since the vector $\va$ is tangent
to ${\M}$, then we have
\begin{equation}\label{Phipna}
\Phi^{\prime \bot}=(f'\va+f\va')^\bot = f\va^{\prime \bot}; \quad
\quad {\Phi^{\prime \bot}}^2=f^2 {\va^{\prime \bot}}^2.
\end{equation}
Taking into account \eqref{apn2} we have ${\va^{\prime\bot}}^2 =
{h'_1}^2-{h'_2}^2$ and therefore ${\Phi^{\prime \bot}}^2=f^2
({h'_1}^2-{h'_2}^2)$. Thus we obtain that the minimal surface
${\M}$ in $\RR^4$ represented by \eqref{W2} is parameterized by
canonical coordinates of the first type if and only if
\begin{equation}\label{can1h}
f^2 ({h'_1}^2-{h'_2}^2)=1
\end{equation}

The last formula implies that the surface ${\M}$ parameterized by
canonical coordinates of the first type has the following
\emph{canonical Weierstass representation}:
\begin{equation}\label{Wcanh}
\Phi: \quad
\begin{array}{rlr}
\phi_1 &=& \ii\ds\frac{\cosh h_1}{\sqrt{{h'_1}^2 - {h'_2}^2}}\\[8mm]
\phi_2 &=&  \ds\frac{\sinh h_1}{\sqrt{{h'_1}^2 - {h'_2}^2}}\\[8mm]
\phi_3 &=&  \ds\frac{\cosh h_2}{\sqrt{{h'_1}^2 - {h'_2}^2}}\\[8mm]
\phi_4 &=& \ii\ds\frac{\sinh h_2}{\sqrt{{h'_1}^2 - {h'_2}^2}}\\
\end{array}
\end{equation}

Conversely, if the pair $(h_1,h_2)$ of holomorphic functions,
determined in a domain in ${\CC}$ satisfy the condition ${h'_1}^2
\neq {h'_2}^2$, then formulas \eqref{Wcanh} give a minimal surface
of general type in $\RR^4$ parameterized by canonical coordinates
of the first type.

If we use the functions $w_1$ and $w_2$ given by \eqref{w1w2},
then the condition \eqref{can1h} gets the form:
\begin{equation}\label{can1w}
f^2 w'_1 w'_2 = 1
\end{equation}

Substituting $h_1$ and $h_2$ in \eqref{Wcanh} by $w_1$ and $w_2$,
respectively, we obtain the following canonical Weierstrass
representation of ${\M}$:
\begin{equation}\label{Wcanw}
\Phi: \quad
\begin{array}{rlr}
\phi_1 &=& \ds\frac{\ii}{\sqrt{w'_1 w'_2}} \cosh \ds\frac{w_1+w_2}{2}\\[6mm]
\phi_2 &=& \ds\frac{1}{\sqrt{w'_1 w'_2}} \sinh \ds\frac{w_1+w_2}{2}\\[6mm]
\phi_3 &=& \ds\frac{1}{\sqrt{w'_1 w'_2}} \cosh \ds\frac{w_1-w_2}{2}\\[6mm]
\phi_4 &=& \ds\frac{\ii}{\sqrt{w'_1 w'_2}} \sinh \ds\frac{w_1-w_2}{2}\\
\end{array}
\end{equation}
Conversely, if $(w_1,w_2)$ is a pair of holomorphic functions
determined in a domain in ${\CC}$, satisfying the condition $w'_1
w'_2 \neq 0$, then formulas \eqref{Wcanw} give a minimal surface
of general type in $\RR^4$ parameterized by canonical coordinates
of the first type.

Finally we obtain a canonical Weierstrass representation of the
type \eqref{W6}. For this aim we use the functions $g_1$ and $g_2$
given by \eqref{g}. After a differentiation of \eqref{g} we get
\begin{equation}\label{gp}
g'_1=e^{w_1}w'_1=g_1 w'_1; \quad g'_2=e^{w_2}w'_2=g_2 w'_2.
\end{equation}
From the above we have
\begin{equation}\label{wp}
w'_1=\frac{g'_1}{g_1}; \quad w'_2=\frac{g'_2}{g_2}.
\end{equation}
Applying  \eqref{ff} in \eqref{can1w}  and \eqref{wp}, we get
$(f2\sqrt{g_1 g_2}\,)^2 \ds\frac{g'_1}{g_1}
\ds\frac{g'_2}{g_2}=1$.

The condition for canonical coordinates of the first type in the
Weierstrass representation \eqref{W6} gets the form:
\begin{equation}\label{can1g}
4f^2 g'_1 g'_2 = 1
\end{equation}

We find $f$ from \eqref{can1g}, replace it into \eqref{W6} and
find the following canonical representation of a minimal surface
of general type:
\begin{equation}\label{Wcang}
\Phi: \quad
\begin{array}{rll}
\phi_1 &=& \ds\frac{\ii}{2}\; \ds\frac {g_1 g_2+1}{\sqrt{g'_1 g'_2}}\\[6mm]
\phi_2 &=& \ds\frac{1}{2}\; \ds\frac {g_1 g_2-1}{\sqrt{g'_1 g'_2}}\\[6mm]
\phi_3 &=& \ds\frac{1}{2}\; \ds\frac {g_1 + g_2}{\sqrt{g'_1 g'_2}}\\[6mm]
\phi_4 &=& \ds\frac{\ii}{2}\; \ds\frac {g_1 - g_2}{\sqrt{g'_1 g'_2}}\\
\end{array}
\end{equation}

Conversely, if $(g_1,g_2)$ is a pair of holomorphic functions
defined in a domain in ${\CC}$ satisfying the condition $g'_1 g'_2
\neq 0$, then formulas \eqref{Wcang} give a minimal surface of
general type in $\RR^4$ parameterized by canonical parameters of
the first type.

\section{Formulas for $K$ and $\varkappa$ in a general Weierstrass representation}

Let ${\M}$ be a minimal surface in $\RR^4$, parameterized by
isothermal coordinates. First we assume that ${\M}$ is given by
the representation \eqref{W5}. In order to obtain formula for $E$,
we use equalities \eqref{EG}, \eqref{adef} and \eqref{a2} and find
\begin{equation}\label{E_fw1w2}
E=|f|^2\cosh(\Re w_1)\cosh(\Re w_2).
\end{equation}
Further we express $\cosh(\Re w_j)$ by means of $g_j$, $j=1,2$ in
view of \eqref{g}:

\begin{equation}\label{cosh_Rew_g}
\cosh(\Re w_j)=\ds\frac{e^{\Re w_j}+e^{-\Re
w_j}}{2}=\ds\frac{e^{2\Re w_j}+1}{2e^{\Re w_j}}=
\ds\frac{|e^{w_j}|^2+1}{2|e^{w_j}|}=\ds\frac{|g_j|^2+1}{2|g_j|}.
\end{equation}
Now making the change \eqref{ff} we get
\begin{equation}\label{E_fg1g2}
E=|f|^2(|g_1|^2+1)(|g_2|^2+1).
\end{equation}

Let us consider the formula \eqref{K_kappa_Phi}. Expressing
$\Phi^{\prime \bot}$ by means of \eqref{Phipna} we get:
$$
K = \ds\frac{-4{\|\Phi^{\prime \bot}\|}^2}{\|\Phi\|^4} =
\ds\frac{-4{\|f\va^{\prime \bot}\|}^2}{\|f\va\|^4} =
\ds\frac{-4{|f|^2\|\va^{\prime \bot}\|}^2}{|f|^4\|\va\|^4} =
\ds\frac{-4{\|\va^{\prime \bot}\|}^2}{|f|^2\|\va\|^4}.
$$
Now using \eqref{mapn2} and \eqref{k1}, we find:
$$
K = \ds\frac{-4(\|\va\|^2\|\va'\|^2-|\bar \va \cdot
\va'|^2)}{|f|^2\|\va\|^6} = \ds\frac{-4k_1}{|f|^2\|\va\|^6}.
$$

In order to obtain a similar formula for $\varkappa$ we use
\eqref{K_kappa_Phi}. We express $\Phi$ by means of $f$ and $\va$
and taking into account \eqref{k2}, we find:
$$
\begin{array}{rll}
\ds \varkappa &=& -\ds\frac{4}{\|\Phi\|^6}\det
(\Phi,\bar\Phi,\Phi^\prime,\overline{\Phi^\prime})
               =  -\ds\frac{4}{\|f\va\|^6}\det (f\va,\bar f \bar \va,f'\va+f\va',\bar f' \bar \va + \bar f \bar {\va'})\\[6mm]
              &=& -\ds\frac{4}{|f|^6\|\va\|^6}\det (f\va,\bar f \bar \va,f\va',\bar f \bar {\va'})
               =  -\ds\frac{4|f|^4}{|f|^6\|\va\|^6}\det (\va,\bar \va,\va',\bar {\va'}) = \ds\frac{4k_2}{|f|^2\|\va\|^6}.
\end{array}
$$
Thus we obtained the following formulas for $K$ and $\varkappa$:
\begin{equation}\label{Kkappa1}
K = \ds\frac{-4k_1}{|f|^2\|\va\|^6}; \quad \varkappa =
\ds\frac{4k_2}{|f|^2\|\va\|^6}.
\end{equation}

Now using \eqref{a2} and \eqref{k12w} we get:
$$
\begin{array}{llr}
K         &=& \ds\frac{-4(|w'_1|^2 \cosh^2 (\Re w_2 ) + |w'_2|^2 \cosh^2 (\Re w_1 ))}{|f|^2 8\cosh^3(\Re w_1)\cosh^3(\Re w_2)}\\[6mm]
\varkappa &=& \ds\frac{ 4(|w'_1|^2 \cosh^2 (\Re w_2 ) - |w'_2|^2
\cosh^2 (\Re w_1 ))}{|f|^2 8\cosh^3(\Re w_1)\cosh^3(\Re w_2)}
\end{array}
$$
From here we find the following formulas for $K$ and $\varkappa$
with respect to the representation \eqref{W5}:
\begin{equation}\label{Kkappa_fw1w2}
\begin{array}{llr}
K         &=& \ds\frac{-1}{2|f|^2\cosh(\Re w_1)\cosh(\Re w_2)}
              \left(\ds\frac{|w'_1|^2}{\cosh^2(\Re w_1)}+\ds\frac{|w'_2|^2}{\cosh^2(\Re w_2)}\right)\\[6mm]
\varkappa &=& \ds\frac{ 1}{2|f|^2\cosh(\Re w_1)\cosh(\Re w_2)}
              \left(\ds\frac{|w'_1|^2}{\cosh^2(\Re w_1)}-\ds\frac{|w'_2|^2}{\cosh^2(\Re w_2)}\right)
\end{array}
\end{equation}

In order to obtain analogous formulas by means of the functions
$g_j$, $j=1;2$, first we note that \eqref{cosh_Rew_g} and
\eqref{wp} imply:
\begin{equation}\label{cosh2_Rew_g}
\ds\frac{|w'_j|^2}{\cosh^2(\Re
w_j)}=\ds\frac{4|g'_j|^2}{(|g_j|^2+1)^2}\quad\quad j=1;2\ .
\end{equation}
Now taking into account the change \eqref{ff} and equality
\eqref{cosh2_Rew_g} we obtain from \eqref{Kkappa_fw1w2} the
following formulas for $K$ and $\varkappa$ with respect to the
representation \eqref{W6}:
\begin{equation}\label{Kkappa_fg1g2}
\begin{array}{llr}
K         &=& \ds\frac{-2}{|f|^2(|g_1|^2+1)(|g_2|^2+1)}
              \left(\ds\frac{|g'_1|^2}{(|g_1|^2+1)^2}+\ds\frac{|g'_2|^2}{(|g_2|^2+1)^2}\right)\\[6mm]
\varkappa &=& \ds\frac{ 2}{|f|^2(|g_1|^2+1)(|g_2|^2+1)}
              \left(\ds\frac{|g'_1|^2}{(|g_1|^2+1)^2}-\ds\frac{|g'_2|^2}{(|g_2|^2+1)^2}\right)\,.
\end{array}
\end{equation}

\section{Formulas for the curvatures $K$, $\varkappa$, $\nu$ and $\mu$ in canonical Weierstrass representation}

Let ${\M}$ be a minimal surface of general type in $\RR^4$,
prameterized by canonical coordinates of the first type. First we
obtain the coefficient $E$ of the first fundamental form in the
canonical Weierstrass representation \eqref{Wcanw}. In the general
form  \eqref{E_fw1w2} we express $f$ under the condition
\eqref{can1w} that the coordinates are canonical of the first type
and find the following formula:
\begin{equation}\label{E_Can_w1w2}
E=\ds\frac{\cosh(\Re w_1)\cosh(\Re w_2)}{|w'_1 w'_2|}.
\end{equation}

In a similar way, we find a formula for $E$ in the case when
${\M}$ is given by the representation \eqref{Wcang}. In view of
\eqref{can1g} we find from the general formula \eqref{E_fg1g2}:
\begin{equation}\label{E_Can_g1g2}
E=\ds\frac{(|g_1|^2+1)(|g_2|^2+1)}{4|g'_1 g'_2|}.
\end{equation}

To obtain formulas for $K$ and $\varkappa$, first let ${\M}$ be
given by means of the representation \eqref{Wcanw}. We find $f$
from the condition \eqref{can1w} and replace it into
\eqref{Kkappa_fw1w2}. Thus we obtain the following formulas for
$K$ and $\varkappa$ in canonical coordinates, with respect to the
representation \eqref{Wcanw}:

\begin{equation}\label{Kkappa_Can_w1w2}
\begin{array}{llr}
K         &=& \ds\frac{- |w'_1 w'_2|}{2\cosh(\Re w_1)\cosh(\Re
w_2)}
              \left(\ds\frac{|w'_1|^2}{\cosh^2(\Re w_1)}+\ds\frac{|w'_2|^2}{\cosh^2(\Re w_2)}\right)\\[6mm]
\varkappa &=& \ds\frac{  |w'_1 w'_2|}{2\cosh(\Re w_1)\cosh(\Re
w_2)}
              \left(\ds\frac{|w'_1|^2}{\cosh^2(\Re w_1)}-\ds\frac{|w'_2|^2}{\cosh^2(\Re w_2)}\right).
\end{array}
\end{equation}

Now let ${\M}$ be given by the representation \eqref{Wcang}. We
find the function $f$ from \eqref{can1g} and replace it into the
general formulas \eqref{Kkappa_fg1g2}. Thus we obtain the
following formulas for $K$ and $\varkappa$ in canonical
coordinates with respect to the representation \eqref{Wcang}:

\begin{equation}\label{Kkappa_Can_g1g2}
\begin{array}{llr}
K         &=& \ds\frac{-8|g'_1 g'_2|}{(|g_1|^2+1)(|g_2|^2+1)}
              \left(\ds\frac{|g'_1|^2}{(|g_1|^2+1)^2}+\ds\frac{|g'_2|^2}{(|g_2|^2+1)^2}\right)\\[6mm]
\varkappa &=& \ds\frac{ 8|g'_1 g'_2|}{(|g_1|^2+1)(|g_2|^2+1)}
              \left(\ds\frac{|g'_1|^2}{(|g_1|^2+1)^2}-\ds\frac{|g'_2|^2}{(|g_2|^2+1)^2}\right).
\end{array}
\end{equation}

Next we find the corresponding formulas for the invariants $\nu$
and $\mu$. Taking into account \eqref{Kkappa_Can_w1w2} we find
\begin{equation}\label{-K+-kappa_Can_w1w2}
\begin{array}{llr}
-K+\varkappa &=& \ds\frac{|w'_1|^3| w'_2|}{\cosh^3(\Re w_1)\cosh(\Re w_2)}\\[6mm]
-K-\varkappa &=& \ds\frac{|w'_1|| w'_2|^3}{\cosh(\Re
w_1)\cosh^3(\Re w_2)}\,.
\end{array}
\end{equation}
Replacing \eqref{-K+-kappa_Can_w1w2} into \eqref{nu_mu_K_kappa},
we obtain formulas for the curvatures $\nu$ and $\mu$ for a
minimal surface of general type given by the representation
\eqref{Wcanw}:

\begin{equation}\label{numu_Can_w1w2}
\begin{array}{llr}
\nu &=& \ds\frac{1}{2}\sqrt{\ds\frac{|w'_1 w'_2|}{\cosh(\Re
w_1)\cosh(\Re w_2)}}
              \left(\ds\frac{|w'_1|}{\cosh(\Re w_1)}+\ds\frac{|w'_2|}{\cosh(\Re w_2)}\right)\\[6mm]
\mu &=& \ds\frac{1}{2}\sqrt{\ds\frac{|w'_1 w'_2|}{\cosh(\Re
w_1)\cosh(\Re w_2)}}
              \left(\ds\frac{|w'_1|}{\cosh(\Re w_1)}-\ds\frac{|w'_2|}{\cosh(\Re w_2)}\right).
\end{array}
\end{equation}
Taking into account \eqref{Kkappa_Can_g1g2}, we get:
\begin{equation}\label{-K+-kappa_Can_g1g2}
\begin{array}{llr}
-K+\varkappa &=& \ds\frac{16|g'_1|^3 |g'_2|}{(|g_1|^2+1)^3 (|g_2|^2+1)}\\[6mm]
-K-\varkappa &=& \ds\frac{16|g'_1| |g'_2|^3}{(|g_1|^2+1)
(|g_2|^2+1)^3}\,.
\end{array}
\end{equation}
Replacing \eqref{-K+-kappa_Can_g1g2} into \eqref{nu_mu_K_kappa} we
obtain formulas for the curvatures $\nu$ and $\mu$ for a minimal
surface of general type given by the representation \eqref{Wcang}:

\begin{equation}\label{numu_Can_g1g2}
\begin{array}{llr}
\nu &=& 2 \ \sqrt{\ds\frac{|g'_1 g'_2|}{(|g_1|^2+1)(|g_2|^2+1)}}
              \left(\ds\frac{|g'_1|}{|g_1|^2+1}+\ds\frac{|g'_2|}{|g_2|^2+1}\right)\\[6mm]
\mu &=& 2 \ \sqrt{\ds\frac{|g'_1 g'_2|}{(|g_1|^2+1)(|g_2|^2+1)}}
              \left(\ds\frac{|g'_1|}{|g_1|^2+1}-\ds\frac{|g'_2|}{|g_2|^2+1}\right).
\end{array}
\end{equation}

With the help of \eqref{Kkappa_Can_g1g2} we can find
transformation formulas for the pair of functions $(g_1,g_2)$
under a motion of the minimal surface ${\M}$ of general type in
$\RR^4$.

Let ${\hat{\M}}$ be another minimal surface of general type in
$\RR^4$, given by the representation \eqref{Wcang} by means of the
pair of functions $(\hat g_1,\hat g_2)$. Both surfaces ${\M}$ and
${\hat{\M}}$ are related by a motion from ${\bf SO}(4,\RR)$ if and
only if they have one and the same curvatures $K$ and $\varkappa$,
calculated with respect to canonical coordinates of te same type.
We note that formulas \eqref{Kkappa_Can_g1g2} coincide with
formulas (2) of \cite{GK1}. Applying Theorem 1 and Theorem 2 in
\cite{GK1} to the curvatures $K$ and $\varkappa$, we obtain that
the surfaces ${\M}$ and ${\hat{\M}}$ are related by a motion from
${\bf SO}(4,\RR)$, if and only if the functions $g_j$ and $\hat
g_j$, \ $j=1,2$ are related by linear fractional transformations
from ${\bf SU}(2,\CC)$:

\begin{equation}\label{hat_g1g2_g1g2}
\hat g_j = \frac{-\bar b_j +\bar a_j\,g_j}{a_j+b_j\,g_j},\quad
a_j=\text{const}, \; b_j=\text{const},\; |a_j|^2+|b_j|^2=1; \;
(j=1;2).
\end{equation}
Replacing $g_j$ with $e^{w_j}$, we obtain transformation formulas
for the pair of functions $(w_1,w_2)$:
\begin{equation}\label{hat_w1w2_w1w2}
e^{\hat w_j} = \frac{-\bar b_j +\bar
a_j\,e^{w_j}}{a_j+b_j\,e^{w_j}},\quad a_j=\text{const}, \;
b_j=\text{const},\; |a_j|^2+|b_j|^2=1; \; (j=1;2).
\end{equation}

\section{Geometric correspondence between minimal surfaces in $\RR^4$, pairs of solutions to the system of natural equations
and pairs of holomorphic functions}

\subsection{Equivalent minimal surfaces in $\RR^4$} In this section we fix a coordinate system
$O(e_1, e_2, e_3, e_4)$ in $\RR^4$, where $\{e_1, e_2, e_3, e_4\}$
is a positive oriented orthonormal quadruple. We suppose that any
minimal surface $(\mathcal M, \text x)$ of general type
$$\mathcal M: \; (u, v) \; \rightarrow \; \text x (u,v); \quad (u, v) \in \mathcal D$$
is defined in a disc $\mathcal D$ with center $(0, 0)$ in $\RR^2
\equiv \CC$ and passes through the point $O$: $\text x(0, 0) = (0,
0, 0, 0)$. The parameters $(u, v)$ are always supposed to be
canonical.

Two minimal surfaces $(\mathcal M, \text x)$ and $(\hat {\mathcal
M}, \hat {\text x})$ of the above type are said to be equivalent
if there exists a disc $\mathcal D_0$, such that
$$\hat {\text x}=A\text x, \quad A \in \mathbf{SO}(4,\RR).$$

We denote by $\mathbf{MS_4}$ the set of equivalence classes of
minimal surfaces of general type in $\RR^4$.

\subsection{Equivalent solutions to the system of natural equations of minimal surfaces of general type in $\RR^4$}

The system of natural equations of minimal surfaces of general
type in $\RR^4$ is the following:
        \begin{equation*}
      \begin{array}{l}
        \ds{(K^2 - \varkappa^2)^{\frac{1}{4}}\, \Delta \ln |\varkappa - K|} = 2(2K - \varkappa)\\
[1mm]
        \ds{(K^2 - \varkappa^2)^{\frac{1}{4}}\, \Delta \ln |\varkappa + K| = 2(2K + \varkappa)}
      \end{array} ; \qquad K<0
    \end{equation*}

Two pairs of solutions $(K,\varkappa )$ and
$(\hat{K},\hat{\varkappa})$ to the above system are said to be
equivalent if there exists a disc $\mathcal D_0$, centered at $(0,
0)$, such that $K= \hat K, \; \varkappa = \hat{\varkappa}$ in
$\mathcal D_0$.

We denote by $\mathbf{SNE_4}$ the set of equivalence classes of
pairs of solutions to the system of natural equations.

\subsection{{Equivalent pairs of holomorphic functions in $\CC$}}
Let $g_k : \mathcal{D}\to \CC $ and $\hat{g_k} :
\mathcal{\hat{D}}\to \CC $,\  $k=1;2$ be two pairs of holomorphic
functions such that $g'_k\neq 0$ и $\hat{g}'_k \neq 0$,\  $k=1,2$.

The two pairs $\{g_1, g_2\}$ and $\{\hat g_1, \hat g_2\}$ are said
to be equivalent if there exists a disc $\mathcal D_0$ such that
        $$\hat g_k = \frac{-\bar b_k +\bar a_k\,g_k}{a_k+b_k\,g_k},\quad
        a_k,b_k=\text{const},\  |a_k|^2+|b_k|^2=1,\quad k=1,2; \;  u+\text{i}v \in \mathcal D_0.$$

We denote by $\mathbf{H^2}$ the set of equivalence classes of
pairs of holomorphic functions.

\subsection{Correspondences between the equivalence classes}

Let $(\mathcal M, \text x)$ be a minimal surface in
$\mathbf{MS_4}$ with Gauss curvature $K$ and normal curvature
$\varkappa$.

Then the correspondence $(\mathcal{M}, \text x)\! \to\!
(K,\varkappa)$ generates a correspondence $\mathbf{MS_4}\! \to\!
\mathbf{SNE_4}$.

This correspondence was obtained by  de Azevero Tribuzy and
Guadalupe \cite{TG}.

Further, let $g_k : \mathcal{D}\to \CC,\  k=1,2 $ be two
holomorphic functions such that $g'_k\neq 0,\  k=1;2$.

Denote by $\Phi$ the vector holomorphic function $\Phi :
\mathcal{D}\to \CC^4 $ defined by the canonical Weierstrass
representation

    \centerline{
    $
  \Phi =
    \left(\ds\frac{\text{i}}{2}\; \ds\frac {g_1 g_2+1}{\sqrt{g'_1 g'_2}}\ ,\
          \ds\frac{1}{2}\; \ds\frac {g_1 g_2-1}{\sqrt{g'_1 g'_2}}\ ,\
          \ds\frac{1}{2}\; \ds\frac {g_1 + g_2}{\sqrt{g'_1 g'_2}}\ ,\
                    \ds\frac{\text{i}}{2}\; \ds\frac {g_1 - g_2}{\sqrt{g'_1 g'_2}} \right)\,.
  $}

Integrating the equality  $\Psi'=\Phi$ we find the function $\Psi
: \mathcal{D}\to \CC^4 $ satisfying the condition $\Psi (0,0) =
(0,0,0,0)$. Then $\text{x}=\Re{\Psi}$ gives a minimal surface
$(\mathcal{M}, \text{x})$ is a minimal surface in $\RR^4$.

Hence the correspondence $(g_1,g_2) \to (\mathcal{M}, \text x)$
generates a correspondence $\mathbf{H^2} \to \mathbf{MS_4}$.

Now, let $g_k : \mathcal{D}\to \CC,\  k=1,2 $ be two holomorphic
functions satisfying the condition $g'_k\neq 0,\  k=1,2$. Then we
find the functions $(K,\varkappa)$ in $\mathcal{D}$ from

\begin{equation*}
   \begin{array}{llr}
   K         &=& \ds\frac{-8|g'_1 g'_2|}{(|g_1|^2+1)(|g_2|^2+1)}
                 \left(\ds\frac{|g'_1|^2}{(|g_1|^2+1)^2}+\ds\frac{|g'_2|^2}{(|g_2|^2+1)^2}\right),\\[4mm]
   \varkappa &=& \ds\frac{ 8|g'_1 g'_2|}{(|g_1|^2+1)(|g_2|^2+1)}
                 \left(\ds\frac{|g'_1|^2}{(|g_1|^2+1)^2}-\ds\frac{|g'_2|^2}{(|g_2|^2+1)^2}\right).
   \end{array}
  \end{equation*}

Thus the correspondence \!$(g_1,g_2)\!\, \to\, (K,\varkappa)$\!
generates a correspondence \! $\mathbf{H^2}\,\to\,\mathbf{SNE_4}$.

This correspondence was obtained by Ganchev and Kanchev in
\cite{GK1}.

Summarizing, we have the following statement:

\begin{thm}
The triangle diagram $(Fig. 1)$ is commutative.
\end{thm}

\begin{figure}[h]\center\epsfysize=5cm\epsffile{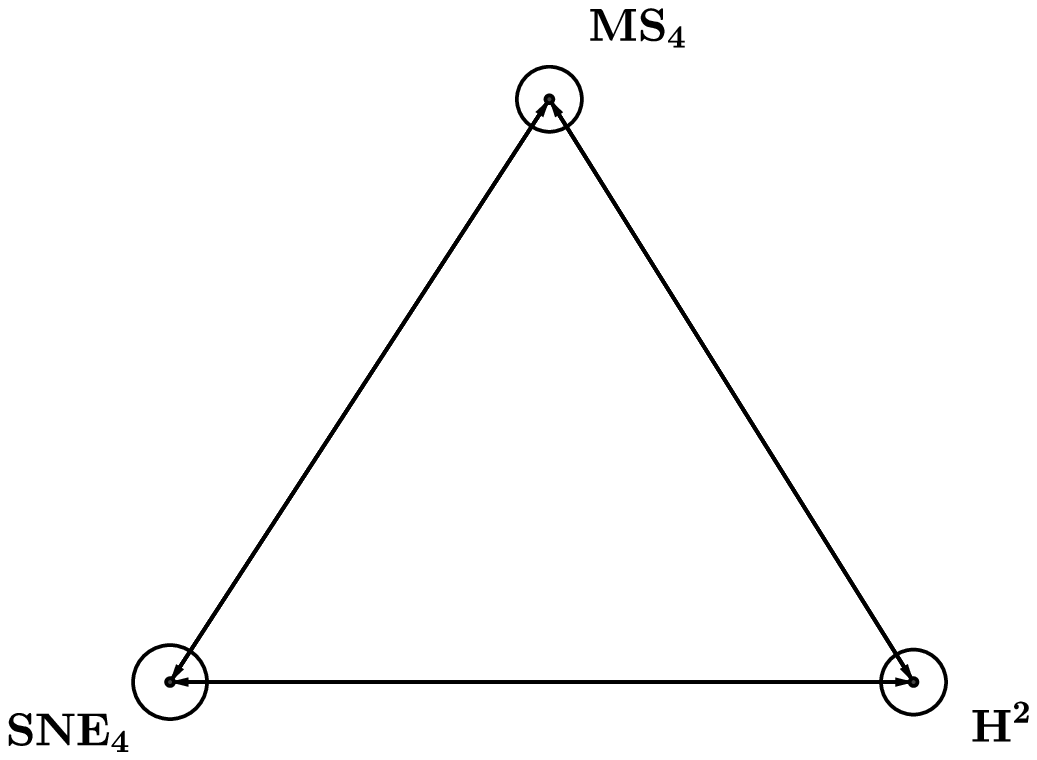}
\center{Fig. 1}
\end{figure}

Finally we shall give a correspondence between minimal surfaces in
$\RR^4$ and pairs of minimal surfaces in $\RR^3$.

First we recall the correspondence between minimal surfaces,
solutions of the natural equation of minimal surfaces and
holomorphic functions.

\section{Geometric correspondence between minimal surfaces in $\RR^3$, solutions to the natural equation
and holomorphic functions}

\subsection{Equivalent minimal surfaces in $\RR^3$.} As in $\RR^4$ we fix a coordinate system
$O(e_1, e_2, e_3)$ in $\RR^3$, where $\{e_1, e_2, e_3\}$ is a
positive oriented orthonormal triple. Let
$$\mathcal M: \; (u, v) \; \rightarrow \; \text x (u,v); \; (u, v) \in \mathcal D$$
be a minimal surface in $\RR^3$ free of umbilical points defined
in a disc $\mathcal D$ with center $(0, 0)$ in $\RR^2 \equiv \CC$.
We consider minimal surfaces passing through the point $O$, so
that $\text x(0, 0) = (0, 0, 0)$. Parameters $(u, v)$ are supposed
to be canonical, i.e. principal and isothermal \cite{G}. If $\nu$
is the positive principal curvature, then the first and the second
fundamental form are given as follows:
$$\mathbf{I} = \frac{1}{\nu}\,(du^2 + dv^2); \qquad \mathbf{II} = du^2 - dv^2.$$
Two minimal surfaces $(\mathcal M, \text x)$ and $(\hat {\mathcal
M}, \hat {\text x})$ of the above type are said to be equivalent
if there exists a disc $\mathcal D_0$ (with center $(0, 0)$), such
that
$$\hat {\text x}=A\text x, \quad A \in \mathbf{SO}(3,\RR).$$

We denote by $\mathbf{MS_3}$ the set of equivalence classes of
minimal surfaces in $\RR^3$.

\subsection{Equivalent solutions to the natural equation of minimal surfaces in $\RR^3$}

The natural equation of minimal surfaces in $\RR^3$ is the
following:
\begin{equation}\label{natural}
\Delta \ln \nu + 2 \nu =0.
\end{equation}

Any solution to the natural equation determines a unique minimal
surface in $\mathbf{MS_3}$.

Two solutions of the natural equation \eqref{natural} are said to
be equivalent if they coincide in a disc $\mathcal D_0$ in $\CC$.

We denote by $\mathbf{SNE_3}$ the set of equivalence classes of
solutions to the natural equation \eqref{natural}.

\subsection{{Equivalent holomorphic functions in $\CC$}}
Let $g : \mathcal{D}\to \CC $ and $\hat{g} : \mathcal{D}\to \CC $,
be two holomorphic functions such that $g'\neq 0$ and $\hat{g}'
\neq 0$. Two holomorphic functions $g$ and $\hat g$ generate one
and the same minimal surface in $\RR^3$ if and only if \cite{OK},
\cite{GK1}:
\begin{equation} \label{ghatg}
\hat g = \frac{-\bar b +\bar a\,g}{a+b\,g},\quad
        a,b=\text{const},\  |a|^2+|b|^2=1,\quad  u+\text{i}v \in \mathcal D_0.
\end{equation}
The two holomorphic functions $g$ and $\hat g$ are said to be
equivalent if they satisfy \eqref{ghatg}.

We denote by $\mathbf{H}$ the set of equivalence classes of
holomorphic functions.

\subsection{Correspondence between the equivalence classes}

Let $ (\mathcal{M}, \text{x})$ be a minimal surfaces in $\RR^3$,
parameterized by canonical parameters. If $\nu$ is the normal
curvature of $\mathcal{M}$, then the correspondence $\mathcal{M}
\to \nu$ generates a correspondence $\mathbf{MS_3} \to
\mathbf{SNE_3}$.

Further, let $g : \mathcal{D}\to \CC $ be a holomorphic function
defined in the disc $\mathcal{D}$ satisfying the condition $g'\neq
0$. Using the canonical Weierstrass representation \cite{G}

\centerline{
    $
  \Phi =
    \left(\ds\frac{{1}}{2}\; \ds\frac {g^2-1}{g'}\ ,\
          \ds-\frac{\text i}{2}\;        \ds\frac {g^2+1}{g'}\ ,\
                                  \ds-\frac {g}{g'}\right)
  $}
\noindent we find the vector holomorphic function $\Psi :
\mathcal{D}\to \CC^3 $ from the equality $\Psi'=\Phi$ and the
condition $\Psi(0, 0)= (0, 0, 0)$. Then $(\mathcal M,\text{x})$,
where $\text{x}=\Re{\Psi}$, is a minimal surface in $\RR^3$.

The correspondence $g \to \mathcal M$ generates a correspondence
$\mathbf{H} \to \mathbf{MS_3}$.

Now let $g : \mathcal{D}\to \CC $ be a holomorphic function
satisfying the condition $g'\neq 0$. This function generates a
solution $\nu: \mathcal D \to \RR $ to the natural equation by
means of the formula \cite{G}
\begin{equation}\label{ni}
\nu =\ds{\frac{4|g^\prime|^2}{(|g|^2+1)^2}}\,.
\end{equation}

The correspondence $g \to \nu$ determines a correspondence
$\mathbf{H} \to \mathbf{SNE_3}$.

Thus we obtained correspondences between $\mathbf {MS_3}, \mathbf
{SNE_3}$ and $\mathbf{H}$:

\begin{figure}[h]\center\epsfysize=5cm\epsffile{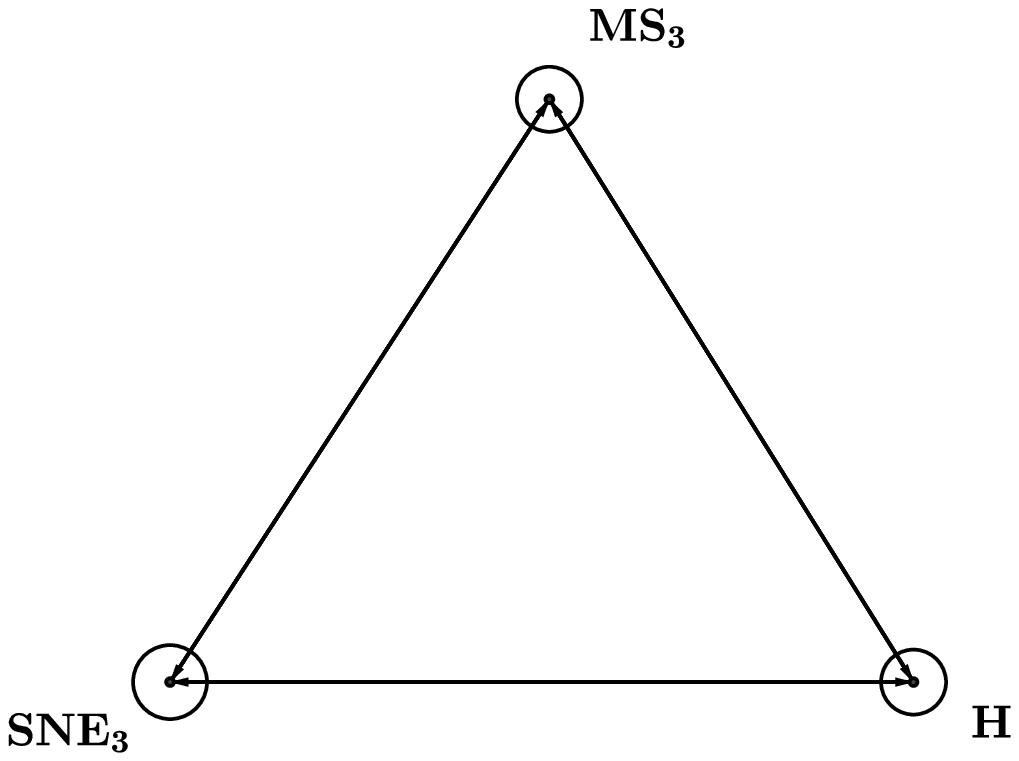}
\center{Fig. 2}
\end{figure}

The triangle diagram (Fig. 2) is commutative and the three
sidelines of the triangle are bijections.

\section{A geometric correspondence between the classes $\mathbf{MS_4}$ and $\mathbf{MS_3} \times \mathbf{MS_3}$}

Let us consider formulas \eqref{Kkappa_Can_g1g2} and \eqref{ni}.
Putting
\begin{equation}
\nu_i =\ds{\frac{4|g_i^\prime|^2}{(|g_i|^2+1)^2}},\,\quad i=1,2.
\end{equation}
we can write the functions $K$ and $\varkappa$ in the form:
\begin{equation*}
      K=-\frac{1}{2} \sqrt{\nu_1\,\nu_2}\,(\nu_1+\nu_2), \quad \varkappa =\frac{1}{2} \sqrt{\nu_1\,\nu_2}\,(\nu_1-\nu_2).
    \end{equation*}
Thus we obtain the statement:
\begin{thm}

\[ \mathbf{SNE_4} \Leftrightarrow \mathbf{SNE_3} \times \mathbf{SNE_3}.\]

\[ \mathbf{MS_4} \Leftrightarrow \mathbf{MS_3} \times \mathbf{MS_3}. \]
\end{thm}

\section{Some applications}
Let us take the holomorphic functions: $g_1=e^{-k_1az}$ and
$g_2=e^{-k_2az}$, where $k_1\neq k_2$ are positive constants,
$a=\cos \alpha + i \sin \alpha$, \; $\alpha =\rm{const}\in [\,0, \pi/4]$ and
$z=u+iv$. Replacing $g_1$ and $g_2$ into \eqref{Wcang} we find a family of
minimal surfaces $\M (k_1,k_2;\alpha)$:
$$
\begin{array}{l}
z_1 = \ds{\frac {1}{k'\sqrt{k_1k_2}}}(\sin 2\alpha \sinh k'p \cos k'q-\cos 2\alpha \cosh k'p \sin k'q),\\
[5mm]
z_2 = \ds{\frac {1}{k'\sqrt{k_1k_2}}}(-\cos 2\alpha \cosh k'p \cos k'q- \sin 2\alpha \sinh k'p \sin k'q),\\
[5mm]
z_3 = \ds{\frac{1}{k''\sqrt{k_1k_2}}}(\cos 2\alpha \sinh k''p \cos k''q+\sin 2\alpha \cosh k''p \sin k''q),\\
[5mm] z_4 = \ds{\frac {1}{k''\sqrt{k_1k_2}}}(-\sin 2\alpha \cosh
k''p \cos k''q+ \cos 2\alpha \sinh k''p \sin k''q),
\end{array}
$$
where $p = u \cos \alpha - v \sin \alpha; \; q = u \sin \alpha + v
\cos \alpha$ and $k'=\frac{k_1+k_2}{2},\; k''=\frac{k_1-k_2}{2}$.

Let us fix $k_1$ and $k_2$. Then we obtain a one-parameter family $\M (\alpha)$.
\begin{itemize}
\item

$\M(0)$ gives the two-parameter family of catenoids in $\RR^4$.

\item
$\M(\pi/4)$ gives the two-parameter family of helicoids in $\RR^4$.

\item
All minimal surfaces $\M(\alpha)$ have the same $K(\alpha)=K(0)$
and $\varkappa(\alpha)=\varkappa(0)$. This implies that any
$\M(\alpha)$ is isometric to $\M(0)$.
\end{itemize}

\begin{rem}
The family $\M(\alpha)$ is the family of the associated with $\M(0)$ minimal surfaces in $\RR^4$.
In some questions in $\RR^4$ the analogue of an isometry  in $\RR^3$ is the notion of a \emph{strong} isometry, i.e.
a deformation of a surface, preserving both $K$ and $\varkappa$.
\end{rem}

\end{document}